\documentclass[12pt]{article}
\usepackage{amsmath,amssymb,amsthm}
\usepackage{amsfonts, amsmath}
\newcommand{\kint}{\smallint} 
\newcommand{\rint}{\int} 
\newtheorem{theorem}{Theorem}[section]
\newtheorem{lemma}{Lemma}[section]
\newtheorem{corollary}{Corollary}[section]

\newtheorem{proposition}{Proposition}[section]
\theoremstyle{definition}
\newtheorem{definition}{Definition}[section]
\theoremstyle{remark}
\newtheorem{remark}{Remarks}[section]
\newtheorem{example}{Example}[section]
\numberwithin{equation}{section}
\begin{document}
\title{\bf Distributions, their primitives and integrals with applications to differential equations}
\author{Seppo Heikkil\"a$^{\dag}$\thanks{Corresponding author\qquad To
appear in {\it Dynamic systems and applications}.} $\ $ and Erik Talvila$^{\ddag}$\\
\small{$^{\dag}$Department of Mathematical Sciences, University of Oulu, Box 3000,} \\
\small{FIN-90014 University of Oulu, Finland}\\
 \small{$^{\ddag}$Department of Mathematics and Statistics,
University of the Fraser Valley,}\\
\small{Abbotsford, BC Canada V2S 7M8}}

\date{2012/09/12}
\maketitle
\title
{}

\begin{abstract} 
In this paper we will  study  integrability of distributions whose primitives are left 
regulated functions and locally or globally integrable in the Henstock--Kurzweil, Lebesgue or Riemann sense. Corresponding spaces of distributions and their primitives are defined and their properties are studied. Basic properties of primitive integrals are derived and applications to systems of first order nonlinear distributional differential equations and to an $m$th order distributional differential equation are presented. The domain of solutions can be unbounded, as shown by concrete examples.
\end{abstract}

\vskip12pt


\baselineskip 16pt

\section{Introduction}\label{S0}

One way of defining an integral is via its primitive. The primitive is a
function whose derivative is in some sense equal to the integrand. For example,
if $f$ and $F$ are functions on a real interval $I$ and $F$ is absolutely continuous,
such that $F'(x) = f(x)$ for almost all $x\in I$, then the Lebesgue integral of $f$
is $\smallint_a^bf(x)\,dx = F(b)- F(a)$  for all $a,\,b\in I$. 
If function $F$ has a pointwise derivative at each point in $I$, 
except for a countable set, then the
derivative is integrable in the Henstock--Kurzweil sense on each compact subinterval of $I$ and
$\smallint_a^b F'(x)\, dx =F(b)-F(a)$ for all $a,\,b\in I$. 
In this sense, the Henstock--Kurzweil integral inverts the pointwise derivative operator. 
There are also  Henstock--Kurzweil integrable functions
for which this fundamental theorem of calculus formula holds and yet these
functions
do not have a pointwise derivative on certain uncountable sets of measure
zero.
A function has a C-integral defined in \cite{BPP00} if
and only if it is everywhere the pointwise derivative of its primitive. 
In this sense, the C-integral is the inverse of the pointwise derivative. It is well-known that the Riemann and Lebesgue
integrals do not have this property. For details, see \cite{ML80}. The
Henstock--Kurzweil integral is equivalent to the Denjoy integral.
 We get the wide Denjoy integral if
 we use the approximate derivative. See, for example, \cite{SS37} for
the definition of the wide Denjoy integral. 

If we use the distributional derivative,
then the primitives need not have any pointwise differentiation properties.
The continuous and regulated primitive integrals  defined in \cite{ET08,ET09} invert the distributional derivatives of  continuous and regulated functions, respectively. See \cite{SH110,SH111,SH12} for applications of these integrals to nonlinear distributional differential equations.

In this paper,
we will study integrability and primitive integrals of distributions on a real interval $I$. 
We say that a distribution $f$ is integrable if $f$ is a distributional derivative of a function, called a primitive of $f$, that is left 
regulated, has a right limit at $\inf I$, and is Henstock--Kurzweil ($HK$) integrable, Lebesgue integrable or Riemann integrable locally on $I$, i.e., on each compact subinterval of $I$. We will show that every
integrable distribution  $f$ also has a left continuous primitive $F:I\to\mathbb R$ that is right continuous at the possible left end point of $I$. Because 
any two such primitives of $f$ differ by a constant, the difference $F(b)-F(a)$ for any two points of $I$ is independent of the particular primitive $F$. This property allows us to define for all $a,\,b\in I$ the primitive integral of $f$ from $a$ to $b$ by 
\begin{equation}\label{E22}
\int_a^bf := F(b)- F(a). 
\end{equation}
%
There is a bijective mapping  $\mathcal F$ between  distributions $f$ and those of their primitives $F$ that have above mentioned one-sided continuity properties, their right limits vanish at $\inf I$, and they are locally integrable in the
$HK$, Lebesgue or Riemann sense. 
In each of these three cases the spaces of primitives have the pointwise partial order $\le$, i.e., $F\le G$ if $F(x)\le G(x)$ for each $x\in I$. The bijection $\mathcal F$ can be used to define a partial order in the corresponding spaces of distributions by $f\preceq g$ if and only if $\mathcal F(f)\le\mathcal F(g)$.
Moreover, if primitives are globally integrable in the $HK$ and Lebesgue cases and $I$ is compact in the Riemann integrable case these spaces can be normed by the Alexiewicz norm $\|\cdot\|_A$ in the $HK$ integrable case, by the $L^1$-norm $\|\cdot\|_1$ in the Lebesgue integrable case, and by the sup-norm $\|\cdot\|_\infty$ in the Riemann integrable case. 
The bijection $\mathcal F$  inherits  norms to the corresponding spaces of distributions by $\|f\|_A=\|\mathcal F(f)\|_A$, $\|f\|_1=\|\mathcal F(f)\|_1$ and $\|f\|_\infty=\|\mathcal F(f)\|_\infty$.     
We will  show that with respect to these partial orderings and norms both the spaces of integrable distributions and their corresponding primitives  form in the
 $HK$ integrable case an ordered normed space, in the Lebesgue integrable case a normed Riesz space, and in the Riemann integrable case a Banach lattice and Banach algebra if $I$ is compact. 
If $I$ is not compact, it can be represented as an increasing denumerable union of compact intervals $I_n$. Thus the spaces of locally integrable primitives can be equipped with the linear metric defined by 
$
d(F_1,F_2)= \underset{n}{\sum}\frac{\|F_1-F_2\|_n}{1+ \|F_1-F_2\|_n},
$
where $\|F\|_n$ denotes the norm of the restriction of $F$ to $I_n$.
  
The Fundamental Theorem of Calculus is valid, so that 
$f=F'$, the primitive derivative of $F=\mathcal F(f)$. This allows us to convert distributional differential equations to integral equations in the spaces of primitive functions.  In \cite{SH122}  this property is applied in the Riemann integrable case  to derive existence results for the unique, smallest, greatest, minimal and/or  maximal solutions of finite systems of first order nonlinear distributional Cauchy problems. Dependence of solutions on the data is also studied, as well as systems of distributional differential equations with impulses and  higher order distributional Cauchy problems.
In section \ref{S8} we generalize to the $HK$ integrable case a uniqueness result and the  existence and comparison results derived in \cite{SH122} for the smallest and greatest solutions of distributional Cauchy systems and  higher order distributional Cauchy problems.
Results of \cite{SH122} dealing with minimal and maximal solutions are extended to the Lebesgue integrable case. Another generalization is that the solution interval can be unbounded, as shown by concrete examples.  

\section{Preliminaries}\label{S1}
\setcounter{equation}{0}

We will first fix some notation for distributions. Let $I$ be a real interval. The space $\mathcal D$ of test
functions are formed by functions of $C^\infty_0(I)$, that is, the smooth functions which, together with all their derivatives, have compact support in  $I$ (cf. \cite{FJ99,TV94}). The support of a function $\phi$ is the closure of the set on which $\phi$ does
not vanish. Denote this as supp($\phi$). 
 There is a notion of continuity in $\mathcal D$.
If $(\phi_n)$ is a sequence in $\mathcal D$, then $\phi_n\to \phi$ in $\mathcal D$ if there is a compact subset $K$ in  $I$ such that for
all $n \in\mathbb N$, supp($\phi_n) \subseteq K$, and for each integer $m\ge  0$, $\phi^{(m)}_n\to \phi^{(m)}$
uniformly
on  $K$ as $n\to\infty$. The distributions on $I$ are the continuous linear functionals
on $\mathcal D$, denoted $\mathcal D'$. If $T \in \mathcal D'$, then $T : \mathcal D \to \mathbb R$ and we write $\left\langle T,\phi\right\rangle\in\mathbb  R$ for
$\phi\in\mathcal D$. If $\phi_n\to \phi$ in $\mathcal D$, then $\left\langle T,\phi_n\right\rangle\to \left\langle T,\phi\right\rangle$ in $\mathbb R$. And, for all $a_1,\, a_2 \in\mathbb R$
and all $\phi_1, \phi_2\in D$, $\left\langle T,a_1\phi_1 + a_2\phi_2\right\rangle = a_1\left\langle T,\phi_1\right\rangle + a_2\left\langle T,\phi_2\right\rangle$.
The differentiation
formula $\left\langle T',\phi\right\rangle = -\left\langle T,\phi'\right\rangle$ ensures that distributions have   derivatives which are distributions. 
Results on distributions can be found in \cite{FJ99}.

A  function $H:I\to\mathbb R$ is left (resp. right) regulated if it has a left (resp. right) limit at each point of $I$ except 
the possible minimum of $I$ (resp. the possible maximum of $I$). Write $H(t-) = \underset{s\to t-}{\lim} H(s)$ and $H(t+) = \underset{s\to t+}{\lim}H(s)$. A function is regulated if it is both right and left regulated. The main difference between regulated functions and right or left regulated functions is that the latter ones may have discontinuities of the second kind, while regulated functions can have only discontinuities of the first kind.  Hence, regulated functions
on a closed interval are bounded while left or right regulated functions
need not be bounded. 
A left regulated function $H$ is left continuous
if $H(t) = H(t-)$ at other points of $I$ than the possible minimum. 
$H$ is said to be countably stepped on a subinterval $[a,b]$ of $I$ if $(a,b]$ is equal to
a countable disjoint union of intervals where $H$ is constant on each interval.

The following lemma presents useful properties for left regulated functions.  

\begin{lemma}\label{L000} Let $H:I\to \mathbb R$ be left regulated. Then 
\newline
(a) $H$ has at most a countable number of discontinuities. 
(b) There is a sequence $(F_n)$ of countably stepped functions on $I$ that $|F_n(t)-H(t)|\le\frac 1n$ for all $n\in\mathbb N$ and $t\in I$. 
(c) $H$ is Lebesgue measurable.
\end{lemma}

{\bf Proof}. (a) Given a compact subinterval $[a,b]$ of $I$ and a positive integer $n$, define $G_n:[a,b]\to [a,b]$ by $G_n(a)=a$,
and for $x\in (a,b]$,
\begin{equation}\label{E402}
G_n(x) = \inf\{y\in[a,x)|\ |H(s)-H(t)|\le \frac 1n\ \hbox{ for all } \ s,\,t\in (y,x)\}, \ x\in(a,b].
\end{equation}
It is easy to verify that $G_n$ is increasing, i.e., $G_n(x)\le G_n(y)$ whenever $a\le x\le y\le b$. Because $H$ is left regulated, then $G_n(x) < x$ for each $x\in(a,b]$.
By \cite[Proposition 1.2.1]{HL94} there is exactly one subset $C_n$ of $[a,b]$ that is inversely well-ordered, i.e., each nonempty subset of $C_n$ has the 
 greatest number, and has the following property:
\begin{equation}\label{E403}
b=\max C_n,\ \hbox{ and $b>x\in C_n$ if and only if } \ x=\inf\{G_n[\{y\in C_n|y> x\}]\}.
\end{equation}
Because $\inf G_n[C_n]$ exists, it is by \cite[Proposition 1.2.1]{HL94} a fixed point of $G_n$.
Since $a$ is the only fixed point of $G_n$, then $\inf G_n[C_n]=a$.
This result and \cite[Theorem 1.2.3]{HL94} imply that $a=\min C_n$. Define
\begin{equation}\label{E413}
D_n :=C_n\setminus\{a\}, \hbox{where $C_n$ is determined by (\ref{E402}), (\ref{E403})}.
\end{equation} 
Because $G_n(x)<x$ for each $x\in D_n$, it follows from dual of \cite[Lemma 1.1.3]{HL94}  that $G_n(x)=\max\{y\in D_n|y < x\}$ for all $x\in D_n$. Thus  
 $(a,b]$ is the disjoint union of half-open intervals $(G_n(x),x]$, $x\in D_n$.   The definition of $G_n$ and the choice of $n$ imply that $|H(s)-H(t)|\le\frac 1n$ for each $s\in(G_n(x),x)$. Thus all the discontinuity points of $H$ in $[a,b]$ belong to the countable set $Z=\bigcup_{m=1}^\infty D_m\cup\{a\}$. This implies the conclusion of (a) because $I$ can be represented as a  denumerable union of its compact subintervals. 

(b) Given a bounded subinterval $(a,b]$ of $I$ and $n\in\mathbb N$, define $F_n:(a,b]\to\mathbb R$ by 
\begin{equation}\label{E102}
F_n(t)=H(x-),\ t\in(G_n(x),x],\ x\in D_n.
\end{equation}
$F_n$ is countably stepped and  $|F_n(t)-H(t)|\le\frac 1n$ for all $t\in (a,b]$. This holds for each $n\in\mathbb N$, so that $(F_n)$ converges to $H$  uniformly on $(a,b]$.     If $\sup I=\infty$ there is $\alpha\in I$ such that
if $x,y\in (\alpha,\infty)$ then $|H(x)-H(y)|<1/n$.  Define $F_n(x)=
\lim_{t\to\infty}H(t)$ for $x\in(\alpha,\infty)$. Now write 
$I\setminus(\alpha,\infty)$ as a disjoint union of intervals $(a,b]$.  Define $F_n$ as above on each such
interval and 
define $F_n(\inf I)=H(\inf I)$.

(c) By (a) the set $Z$ of discontinuity points of $H$ is a null set, whence $H$ is Lebesgue measurable.
\qed

Applying results of Lemma \ref{L000} we obtain the following integrability criteria for a left regulated function.

\begin{lemma}\label{L001} Let $H:I\to \mathbb R$ be left regulated. Then 
(a) $H$ is locally Riemann integrable if and only if $H$ is locally bounded.  
(b) $H$ is locally Lebesgue integrable if and only if for each subinterval $[a,b]$ of $I$ the function $F_n$, defined by (\ref{E102}), is Lebesgue integrable for some $n\in\mathbb N$.
(c) $H$ is locally $HK$ integrable if and only if for each subinterval $[a,b]$ of $I$ the function $F_n$, defined by (\ref{E102}), is $HK$ integrable for some $n\in\mathbb N$.
\end{lemma}

{\bf Proof}. (a) Because $H$ by Lemma \ref{L000} is continuous almost everywhere on $I$, then $H$ is locally Riemann integrable if and only if $H$ is locally bounded (see \cite{M59}). 

(b) and (c) Let $[a,b]$ be a subinterval of $I$. By Lemma \ref{L000} (b) the function $F_n$ defined by (\ref{E102}) is countably stepped, and hence Lebesgue measurable.
If $H$ is locally $HK$ integrable, then $F_n$ is bounded above and below by $HK$ integrable functions $H\pm \frac 1n$, whence $F_n$ is $HK$ integrable  (see \cite[Theorem 2.5.16]{LPY00}). Conversely, if $F_n$ is $HK$ integrable, then $H$ is bounded above and below on $[a,b]$ by $HK$ integrable functions $F_n\pm \frac 1n$. Because $H$ is also Lebesgue measurable by Lemma~\ref{L001} (c), then $H$ is $HK$ integrable on $[a,b]$. 
The above reasoning holds also when $HK$ integrability is replaced by Lebesgue integrability.    
\qed  

The following lemma, which is a consequence, e.g., of \cite[Lemma 1.12]{CH11}, presents a sufficient condition for local $HK$ integrability.

\begin{lemma}\label{L002} Let $I$ be an interval in $\mathbb R$. 
Given a function $G:I\to \mathbb R$, suppose that there exists a continuous
function $F: I\to \mathbb R$ and  a countable subset $Z$ of $I$ such
that $F$ is differentiable in  $I\setminus Z$, and $F'(t)=G(t)$ for all $t\in I\setminus Z$. Then $G$ is locally $HK$
integrable on $I$, and $\int_a^bG(t)\,dt = F(b)-F(a)$ for all $a,\,b\in I$.
\end{lemma}

The next result follows from  \cite[(8.6.4)]{Die60}. 

\begin{lemma}\label{L003} Let $(G_m)_{m=1}^\infty$ be a sequence of functions from an interval $I\subseteq \mathbb R$ into $\mathbb R$. Suppose that, for each $m\in\mathbb N$, there exists a continuous function $F_m:I\to\mathbb R$ and a countable subset $Z_m$ of $I$ such that  $F_m$ is differentiable in  $I\setminus Z_m$, and $F_m'(t)=G_m(t)$ for all $t\in I\setminus Z_m$. Suppose in addition that
\begin{enumerate}
\item[(i)] there is a point $t_0\in I$ such that $(F_m(t_0))$ converges in $\mathbb R$;
\item[(ii)] for every point $t\in I$ there is a neighbourhood $B(t)$ with respect to $I$ such that in $B(t)$ the sequence $(G_m)$ converges uniformly. 
\end{enumerate}
Then for each $t\in I$, the sequence $(F_m)$ converges uniformly in $B(t)$; and if we put $G(t)=\underset{m\to\infty}{\lim}G_m(t)$ and $F(t)=\underset{m\to\infty}{\lim}F_m(t)$, then $F'(t)=G(t)$ at every point $t$ of $I\setminus \underset{n}{\cup}Z_m$.
\end{lemma}

Lemmas \ref{L002} and \ref{L003} are used in Example \ref{Ex601} to verify the $HK$ integrability of a left regulated function that has a discontinuity of the second kind at every rational point.

The definition of integrability and the primitive integral in the $HK$ integrable case is based on the following result, where $\smallint$ denotes the Henstock--Kurzweil integral.

\begin{lemma}\label{L00} Suppose that $H:I\to\mathbb R$ is left regulated, that $\phi_n\to\phi$ in $\mathcal D$, that $K$ is the compact subset of $I$ as in the definition of $\phi_n\to\phi$, and that  $K\subseteq [a,b]\subseteq I$. Then
$\kint_a^b H(t)\phi_n(t)\,dt\to \kint_a^b H(t)\phi(t)\,dt$ if and only if $H$ is $HK$ integrable on $[a,b]$.
\end{lemma}

{\bf Proof}. Since $H$ is left regulated, it is  Lebesgue measurable by Lemma \ref{L000}. Since the sequence $(\phi_n-\phi)$ converges to $0$ in $\mathcal D$, the sequence $(\phi'_n-\phi')$ converges to $0$ uniformly on $[a,b]$. Thus the sequence $(\phi_n)$ is of uniform bounded variation, and converges uniformly to $\phi$.
The conclusion follows from \cite[Corollary 3.2]{ET99}, since the hypotheses of it  are valid by the above proof when $f=H$, $(g_n)=(\phi_n)$ and $g=\phi$.
\qed

\section {The $LD$ primitive integral and basic properties}\label{S2}
\setcounter{equation}{0}

In this section we will study integrability and the integral of distributions on a real interval $I$ having locally or globally $HK$ integrable primitives.

\subsection{$LDP$ integrability and the $LD$ primitive integral}\label{S21}
 
We will first describe the spaces of primitives for the $LD$ primitive
integral on  $I$. 
Denote 
\begin{equation}\label{E100}
\begin{aligned}
\mathcal D^{lr}(I) =&\{H : I\to \mathbb R | \hbox{$H$ is left regulated, locally $HK$ integrable,}\\ 
&\hbox{and has a right limit at $\inf I$}\},\\
\mathcal D^{lc}(I) =&\{G\in \mathcal D^{lr}(I) | \hbox{$G$ is left continuous, and $G(\min I) = G(\min I+)$ if $\min I$ exists}\},\\
\mathcal D^{lc}_0(I)=&\{F\in \mathcal D^{lc}(I)| F(\inf I+)=0\}. 
\end{aligned}
\end{equation}
The same notations are also used when local $HK$ integrability is replaced by global integrability. 

Let $H \in \mathcal D^{lr}(I)$. We will prove in  Theorem \ref{T0} that $H$ uniquely 
determines a distribution, also denoted by  $H$, on $I$ by   
\begin{equation}\label{E0011}
\langle H,\phi\rangle = \int_a^bH(t)\phi(t)\,dt, \quad \phi\in\mathcal D, \ \hbox{where supp($\phi)\subseteq [a,b]\subseteq I$}, 
\end{equation} 
where $\smallint$ denotes the Henstock--Kurzweil integral. Define functions $G\in \mathcal D^{lc}(I)$ and $F\in\mathcal D^{lc}_0(I)$ by
\begin{equation}\label{E011} 
G(x)= \begin{cases} H(x-), \ x>\inf I,\\ H(x+), \ x=\min I \hbox{ if $\min I$ exists},\end{cases} \quad F(x)=G(x)-G(\inf I+), \ x\in I.
\end{equation}
Replacing $H$ in (\ref{E0011}) by the so defined functions $G$ and $F$ we get distributions $G$ and $F$.
Because $H$ has by  Lemma \ref{L000} only a countable number of discontinuity points, then 
all the distributions $F$, $G$ and $H$ are equal.
All these three distributions have the same derivative which is itself a distribution. This is known as the distributional
derivative or weak derivative. 
We will usually denote the distributional
derivative of a distribution $F$ by $F'$ and the possible pointwise derivative of $F$ by $F'(t)$.

As a consequence of Lemma \ref{L00} we obtain

\begin{theorem}\label{T0} Every function of $\mathcal D^{lr}(I)$ determines  a unique distribution on $I$ by (\ref{E0011}). 
\end{theorem}

{\bf Proof}. Let $H\in \mathcal D^{lr}(I)$. Since $H$ is locally $HK$ integrable and left regulated, then (\ref{E0011}) defines
by Lemma \ref{L00} a continuous function $H:\mathcal D\to\mathbb R$ by 
$\phi\mapsto \kint_a^bH(s)\phi(t)\,dt$, where supp($\phi)\subseteq[a,b]\subseteq I$.
$H$ is linear because of linearity properties of Henstock--Kurzweil integrals. Thus $H\in \mathcal D'$. 
\qed

A distribution $f$ on $I$ is called $LDP$ integrable on $I$ 
if it is the distributional derivative of some primitive $H \in \mathcal D^{lr}(I)$, that is,
for all $\phi\in\mathcal D$ we have 
\begin{equation}\label{E20}
\left\langle f,\phi\right\rangle = \left\langle H',\phi\right\rangle = -\left\langle H,\phi'\right\rangle = -
\int_a^b H(t)\phi'(t)\, dt, \ \hbox{where supp($\phi')\subseteq[a,b]\subseteq I$.}
\end{equation}
 The last integral is a Henstock--Kurzweil integral. 
 Denote
\begin{equation}\label{E21}
\mathcal A_D(I) = \{f \in \mathcal D' | f \ \hbox{ is $LDP$ integrable on } \ I\}.
\end{equation} 
 Linearity of the distributional derivative shows that $\mathcal A_D(I)$ is a linear subspace
of $\mathcal D'$.
If $f$ is $LDP$ integrable with a primitive $H\in \mathcal D^{lr}(I)$, then (\ref{E011})  determines a primitive $F$ of $f$ in $\mathcal D^{lc}_0(I)$. 

If $f\in \mathcal A_D(I)$, $F$ is a primitive of $f$ in $\mathcal D^{lc}(I)$  and $a,\,b\in I$, define the $LD$ primitive integral of $f$ from $a$ to $b$ by (\ref{E22}). 

\subsection{Basic properties of the $LD$ primitive integral}\label{S22}

We now present some of the basic properties of the $LD$ primitive integral. Linear combinations are defined by 
$\left\langle a_1f_1 +
a_2f_2,\phi\right\rangle = \left\langle a_1F'_1 + a_2F'_2 ,\phi\right\rangle$ for $\phi \in \mathcal D$; $a_1, a_2 \in\mathbb R$; $f_1, f_2 \in \mathcal A_D(I)$ with primitives
$F_1,F_2 \in\mathcal D^{lc}_0(I)$.

\begin{theorem}\label{T2} (Basic properties of the integral). (a) The $LD$ primitive integral is unique.
(b) Addivity over intervals. If $f \in \mathcal A_D(I)$, then for all $a\le c\le b$ we
have $\int_a^c f + \int_c^b f =\int_a^b f$. 
(c) Linearity. If $f_1, f_2 \in  \mathcal A_D(I)$ and $a_1, a_2 \in\mathbb R$, then $a_1f_1 + a_2f_2 \in \mathcal A_D(I)$ and $\int_a^b(a_1f_1 + a_2f_2) = a_1\int_a^bf_1 + a_2\int_a^bf_2, \quad a,\,b\in I$.
(d) Reverse limits of integration.
Let $a,\,b\in I$. Then
$\int_b^a f = -\int_a^b f$.
\end{theorem}

{\bf Proof}. (a) To prove the $LD$ primitive integral is unique we need to prove primitives
in $\mathcal D^{lc}_0(I)$
 are unique. Suppose $F,G \in\mathcal D^{lc}_0(I)$ and $F' = G'$. Then $(F - G)' = 0$ and
the only solutions of this distributional differential equation are the constant
distributions \cite[Section 2.4]{FJ99}. The only constant distribution in $\mathcal D^{lc}_0(I)$ is the
zero function.

(b) Note that $[F(b)- F(c)] + [F(c) - F(a)] = F(b) - F(a)$.
 
(c) Since $a_1f_1+a_2f_2 = (a_1F_1+a_2F_2)'$, where $F_i \in\mathcal D^{lc}_0(I)$, $F_i'=f_i$, we have
$$
\begin{aligned}&\int_a^b(a_1f_1+a_2f_2) = (a_1F_1+a_2F_2)(b)-(a_1F_1+a_2F_2)(a)\\ 
&= a_1F_1(b) + a_2F_2(b)-a_1F_1(a) - a_2F_2(a)
= a_1\int_a^bf_1 + a_2\int_a^bf_2.
\end{aligned}
$$

(d) $\int_b^a f = F(a) - F(b) = -[F(b) - F(a)] = -\int_a^bf$. 
\qed

\begin{theorem}\label{T3} (Fundamental theorem of calculus). 
Let $f\in \mathcal A_D(I)$, $G\in \mathcal D^{lc}(I)$, $a\in I$ and
 $c\in\mathbb R$. Then  $G'=f$ and $G(a)=c$ if and only if  $G(x)=c + \int_{a}^x f$ for every $x\in I$.
\end{theorem}

{\bf Proof}. 
If $G'=f$ and $G(a)=c$, it follows from (\ref{E011}) and (\ref{E22}) that 
$$
\begin{aligned}
G(x)&=c+G(x)-G(a)=c+(G(x)-G(\inf I+))-(G(a)-G(\inf I+))\\
&=c+F(x)-F(a)=c+\int_a^xf.
\end{aligned}
$$
Conversely, let  $G(x)=c+\int_a^xf$. Then $G(a)=c+\int_a^af=c$ by (\ref{E22}). Since $F_a(x)=\int_a^xf$ is a primitive of
$f$ in $\mathcal D^{lc}(I)$, then $F_a'=f$. Thus $G'=(c+F_a)'=0+F_a'=f$.
\qed

The result of Theorem \ref{T3} can be used to convert distributional initial value problems
into integral equations.  

As a consequence of the definitions of  Theorem \ref{T2} and the definitions of $\mathcal A_D(I)$ and $\mathcal D^{lc}_0(I)$ we obtain

\begin{corollary}\label{C2} The  mapping $\mathcal F$, defined by
\begin{equation}\label{E123} 
\mathcal F(f)=F, \ f\in\mathcal A_D(I),\ \hbox{where $F$ is the primitive of $f$ in}\ \mathcal D^{lc}_0(I),
\end{equation}
is a linear isomorphism from $\mathcal A_D(I)$ to $\mathcal D^{lc}_0(I)$.
\end{corollary}

As with the Henstock--Kurzweil integral, there are no improper integrals.

\begin{theorem}\label{T4} (Hake theorem). Suppose $f \in \mathcal D'$, and that $f \in\mathcal A_D([x,y])$ for every proper subinterval $[x,y]$ of $[a,b]\subseteq I$. If for some $c \in(a,b)$,
$\underset{x\to a+}{\lim}
\int_x^c f$ exists, then $f\in\mathcal A_D[a,c]$, and $\int_a^c f=\underset{x\to a+}{\lim}
\int_x^c f$, and if    $\underset{y\to b-}{\lim}\int_c^y f$ exists, then  $f\in\mathcal A_D[c,b]$, and $\int_c^bf=\underset{y\to b-}{\lim}\int_c^y f$.
\end{theorem}

\subsection{Order and norm properties}\label{S23}

In $\mathcal D^{lc}_0(I)$, there is the partial order: $F \le G$ if and only if $F(x) \le G(x)$ for all
$x \in I$. 
Since the mapping $\mathcal F$, defined by (\ref{E123}) is a linear isomorphism from $\mathcal A_D(I)$ to $\mathcal D^{lc}_0(I)$, we can define a partial order $\preceq$ in $\mathcal A_D(I)$ as follows.
For $f,\, g\in\mathcal A_D(I)$, define  $f \preceq g$ if $\mathcal F(f)\le \mathcal F(g)$ in $\mathcal D^{lc}_0(I)$. In particular, if $\min I$ exists, then 
\begin{equation}\label{E30}
f \preceq g\ \hbox{ if and only if $\int_{\min I}^xf\le \int_{\min I}^xg$ for each } \ x\in I.
\end{equation}
Thus $f\preceq g$ if and only if $F \le G$, where $F$ and $G$ are the
respective primitives in $\mathcal D^{lc}_0(I)$. 
If $\preceq$ is a binary operation on set $E$, then
it is a partial order if for all $x,\, y,\, z \in E$ it is reflexive ($x \preceq x$), antisymmetric
($x \preceq y$ and $y \preceq x$ imply $x = y$) and transitive ($x \preceq y$ and $y \preceq z$ imply $x \preceq z$).
If $E$ is a vector space and $\preceq$ is a partial order on $E$, then $E$
is an ordered vector space if for all $x,\, y,\, z \in S$

(1) $x \preceq y$ implies $x + z \preceq y + z$.

(2) $x \preceq y$ implies $kx \preceq ky$ for all $k \in\mathbb R$ with $k \ge 0$.

If $x \preceq y$, we write $y \succeq x$. 

\begin{theorem}\label{T5} (Ordered vector space). (a) Both  $\mathcal D^{lc}(I)$ and $\mathcal D^{lc}_0(I)$ are ordered vector spaces. (b)  $\mathcal A_D(I)$ is  order isomorphic to $\mathcal D^{lc}_0(I)$. 
\end{theorem}

{\bf Proof}. (a) 
The following properties follow immediately from the definition. If $F \le G$ in $\mathcal D^{lc}(I)$, 
then for all $H \in\mathcal D^{lc}(I)$ we have $F + H \le G + H$. If $F \le G$ and $k\geq 0$ then
$kF \le kG$.  Hence, $\mathcal D^{lc}(I)$ is an ordered vector space, and so is $\mathcal D^{lc}_0(I)$ as a linear subspace of $\mathcal D^{lc}(I)$.

(b) If $f,\, g \in\mathcal A_D(I)$ and $f \preceq g$, then $\mathcal F(f) \le \mathcal F(g)$. Let $h \in\mathcal A_D(I)$. Then, $\mathcal F(f) + \mathcal F(h) \le \mathcal F(g) + \mathcal F(h)$.
But then $(\mathcal F(f) + \mathcal F(h))' = \mathcal F(f)' + \mathcal F(h)' = f + h \preceq g + h$. If $k \in\mathbb R$ and $k \geq 0$, then
$(k\mathcal F(f))' = k\mathcal F(f)' = kf$ so $kf \preceq kg$. 
Then $\mathcal A_D(I)$ is an ordered vector space that is order isomorphic to $\mathcal D^{lc}_0(I)$.
\qed

A vector space $E$ equipped with a partial order $\preceq$ and a norm $\|\cdot\|$ is said to be an ordered normed space if the order cone $E_+=\{x\in E|x\succeq 0\}$ is a closed subset of $E$ in its norm topology.  

Assume next that the functions of $\mathcal D^{lc}(I)$ and $\mathcal D^{lc}_0(I)$ are $HK$ integrable. We will show that $\mathcal D^{lc}(I)$ and $\mathcal D^{lc}_0(I)$, ordered pointwise, are ordered normed spaces with respect to the Alexiewicz norm:  
$$
\|F\|_A=\underset{[a,b]\subseteq I}{\sup}\left |\kint_a^bF(t)\,dt\right |, \quad F\in\mathcal D^{lc}(I).
$$
Using the isomorphism $\mathcal F$ we define a norm in $\mathcal A_D(I)$ by $\|f\|_A=\|\mathcal F(f)\|_A$. Equivalently,  
\begin{equation}\label{E301}
\|f\|_A= \underset{[a,b]\subseteq I}{\sup}\left |\kint_a^bF(t)\,dt\right |, \quad  \hbox{ where }\ F=\mathcal F(f).
\end{equation}

\begin{theorem}\label{T6} (Ordered normed space). Assume that the functions of $\mathcal D^{lc}(I)$ and $\mathcal D^{lc}_0(I)$ are $HK$ integrable. (a) $\mathcal D^{lc}(I)$ and $\mathcal D^{lc}_0(I)$ are ordered normed spaces
with respect to the Alexiewicz norm and pointwise order. (b) $\mathcal A_D(I)$ is an ordered normed space with respect to partial order and norm defined by (\ref{E30}) and (\ref{E301}). (c) $\mathcal A_D(I)$ and $\mathcal D^{lc}_0(I)$ are isometrically 
isomorphic.  The integral provides a linear isometry.
\end{theorem}

{\bf Proof}. (a) To prove that the Alexiewicz norm is a norm in $\mathcal D^{lc}(I)$ and in $\mathcal D^{lc}_0(I)$, 
 first note they are linear subspaces of the Denjoy space $D(I)$ of all $HK$ integrable functions from $I$ to $\mathbb R$.
And, if $F\in \mathcal D^{lc}(I)$ such that $\|F\|_A =0$, then 
then $F(x) = 0$ for almost all $x \in I$. But $F$ is left continuous in $I\setminus\{\inf I\}$ and right continuous at the possible  $\min I$. So if there were
$b \in I$ such that $F(b) \ne 0$ then there is an interval $(a,b]$ if $b>\min I$ or an interval $[b,c)$ if $b=\min I$ in which $F$ does not vanish,
which is a contradiction. Thus $F(x) = 0$ for all $x \in I$. Positivity, homogeneity
and the triangle inequality are inherited from $D(I)$.

(b) Because  $HK$ integrability and strong $HK$ integrability (called HL integrability in \cite{CH11})  are equivalent for real functions by \cite[Proposition 3.6.6]{Sye05}, it follows from \cite[Lemma 9.29]{CH11} that the cone $D(I)_+=\{G\in D(I)| G(x)\ge 0$ for a.e. $x\in I\}$ is closed with respect to the Alexiewicz norm. Hence, if $\|F_n-F\|_A\to 0$ in $\mathcal D^{lc}(I)$ and $F_n\in \mathcal D^{lc}(I)_+=\{G\in\mathcal D^{lc}(I)|G(x)\ge 0$ for all $x\in I\}$, then $F(x)\ge 0$ for a.e. $x\in I$. Because $F$ is left continuous in $I\setminus\{\min I\}$ and right continuous at $\min I$  one can show (cf. the proof of (a)) 
that 
$F$ cannot have negative values, whence $F\in \mathcal D^{lc}(I)_+$. This proves that $\mathcal D^{lc}(I)_+$ is closed in
the Alexiewicz norm topology of $\mathcal D^{lc}(I)$, so $\mathcal D^{lc}(I)$ is an ordered normed space. 
The proof that $\mathcal D^{lc}_0(I)$ is an ordered normed space is similar.  

(c) The conclusion is a direct consequence of (a) because $\mathcal F$ is an order isomorphism from $\mathcal A_D(I)$ to $\mathcal D^{lc}_0(I)$,  and  $\|f\|_A=\|\mathcal F(f)\|_A$ for all $f\in\mathcal A_D(I)$.
\qed

\subsection{Integration by parts, dual space}\label{S24}
It is well-known that if $f\in D([a,b])$ then the pointwise product
$fg$ is also Henstock--Kurzweil integrable if $g$ is of bounded variation.
Let $F(x)=\int_a^xf$.
The integration by parts formula is then defined in terms of
a Riemann--Stieltjes integral via
$\int_a^b f(x)g(x)\,dx=F(b)g(b)-\int_a^bF(x)\,dg(x)$.
See \cite{ML80}. 
The functions
of (essential) bounded variation also form the dual space of $D([a,b])$.  
We will
see analogues of these results for the $LD$ primitive integral.

Let $\mathcal{BV}([a,b])$ be the functions of bounded variation on
$[a,b]$.
\begin{definition}\label{D24.0}
Let $c\in[a,b]$.  Define 
\begin{equation}
\mathcal{IBV}_c([a,b])=\{g:[a,b]\to\mathbb R\mid
g(x)=\int_c^xh(t)\,dt\quad\text{ for some }h\in
\mathcal{BV}([a,b])\cap L^1([a,b])\}.
\end{equation}
\end{definition}
The $L^1$ condition is redundant if $[a,b]$ is a compact interval.
Note that functions in $\mathcal{IBV}_c([a,b])$ vanish at $c$,
are Lipschitz continuous on $[a,b]$ and are hence in $\mathcal{BV}([a,b])$.

\begin{definition}\label{D24.1}
Let $f\in\mathcal A_D([a,b])$ with primitive $F\in\mathcal D^{lc}_0([a,b])$.
Let $g\in\mathcal{IBV}_c([a,b])$ for some $c\in[a,b]$.
Define the integration by parts formula
$\int_a^bfg=F(b)g(b)-\int_a^bF(t)g'(t)\,dt$.
\end{definition}
There is no way of proving the integration by parts definition, although it
clearly
holds if $f\in D([a,b])$.  However,
we can use a sequential approach to justify it since the $C^1$ functions
are dense in $D([a,b])$.  Functions in $\mathcal{IBV}_a$ and $\mathcal{IBV}_c$
differ by a constant so we just need consider $\mathcal{IBV}_a$.
\begin{theorem}\label{theorem24.2}
Let $f\in\mathcal A_D([a,b])$ with primitive $F\in\mathcal D^{lc}_0([a,b])$.
Let $g\in\mathcal{IBV}_a([a,b])$. (a)  Then
$|\int_a^bfg|\leq|F(b)||g(b)|+\|F\|_A(\|g'\|_\infty+Vg')$.
(b) Suppose $(F_n)\subset C([a,b])\cap C^1((a,b))$ such that $F_n(a)=0$,
$F_n(b)=F(b)$ and 
$\|F'_n-f\|_A\to 0$.
Let $H(x)=
F(x)g(x)-\int_a^xF(t)g'(t)\,dt$ and let 
$H_n(x)=\int_a^xF'_n(t)g(t)\,dt$.
Then
$|H(b)-H_n(b)|\leq
\|F-F_n\|_AVg'\to 0$ 
as $n\to\infty$.
\end{theorem}

{\bf Proof}. (a) The inequality
$\left|\int_a^b F\phi\right|\leq \|F\|_A(\inf|\phi|+V\phi)$ holds for all
$F\in D([a,b])$ and $\phi\in\mathcal{BV}([a,b])$.  See \cite{CD89} and 
\cite[Lemma~24]{ET02}.

(b) By the multiplier
result above, the product $Fg\in\mathcal D^{lc}_0([a,b])$.  Since $g'$ is almost everywhere
equal to a function of bounded variation we also have $Fg'\in D([a,b])$.
The function $x\mapsto \int_a^xF(t)g'(t)\,dt$ is continuous on $[a,b]$ 
and vanishes at $a$.
Hence, $H$ exists on $[a,b]$.  There exist functions $F_n$ satisfying the endpoint conditions
since $F$ has limits at $a+$ and $b-$.  Each function $H_n$ is continuous
on $[a,b]$.
The inequality follows from (a).
\qed

For compact intervals there is convergence in the Alexiewicz norm to the
integration by parts formula.
\begin{theorem}\label{theorem24.3}
Let $[a,b]$ be a compact interval.
Let $f\in\mathcal A_D([a,b])$ with primitive $F\in\mathcal D^{lc}_0([a,b])$.
Let $g\in\mathcal{IBV}_a([a,b])$.
Suppose $(F_n)\subset C([a,b])\cap C^1([a,b])$ such that
$\|F'_n-f\|_A\to 0$.
Let $H(x)=
F(x)g(x)-\int_a^xF(t)g'(t)\,dt$ and let 
$H_n(x)=\int_a^xF'_n(t)g(t)\,dt$.
Then $H,H_n\in \mathcal D^{lc}_0([a,b])$ and
$\|H-H_n\|_A\to 0$.
\end{theorem}

{\bf Proof}.
The proof of Theorem~\ref{theorem24.2}(b) shows that
$H,H_n\in \mathcal D^{lc}_0([a,b])$.
Now, 
\begin{equation}
\|H'-H'_n\|_A=\|H-H_n\|_A\leq 2\sup_{\beta\in[a,b]}\left|
\int_a^\beta\left[H(x)-H_n(x)\right]dx\right|.
\end{equation}
Let $\beta\in[a,b]$.
Use (a) of Theorem~\ref{theorem24.2}.  Then
\begin{equation}
\begin{aligned}
\left|
\int_a^\beta\left[H(x)-H_n(x)\right]dx\right|&\leq
\left|\int_a^\beta\left[F(x)-F_n(x)\right]g(x)\,dx\right|\\
&\qquad+\left|
\int_a^\beta\left[F(t)-F_n(t)\right]g'(t)\int_a^\beta\chi_{[t,\beta]}(x)\,dx\,dt
\right|\\
&\leq \|F-F_n\|_A\left[(\|g\|_\infty+Vg)+(b-a)(\|g'\|_\infty+Vg')\right].
\end{aligned}
\end{equation}
The order of the iterated integrals can be interchanged by 
\cite[Theorem~57]{CD89}.  It now follows that $\|H-H_n\|_A\to0$.
\qed

Notice that integration by parts defines a product
$\mathcal A_D([a,b])\times\mathcal{IBV}_a([a,b])\to\mathcal A_D([a,b])$
given by 
$fg=H'$.  A similar type of definition was used
in \cite{ET09}.  Various properties of this product (or bimodule) were
proved in Theorem~18 of that paper.  Since $\mathcal{IBV}_a([a,b])$ 
is closed under pointwise products,
the same results hold here with
similar proofs.

Note that if $J$ is a subinterval of $[a,b]$ then $\int_Jf$ is not
defined since $\chi_J\not\in\mathcal{IBV}_c$.

Now we consider the dual space of $\mathcal A_D([a,b])$.
First, denote the functions of essential bounded variation by 
$\mathcal{EBV}([a,b])$.
A function $f\in\mathcal{EBV}([a,b])$ if $f=g$ a.e. for some $g\in
\mathcal{BV}([a,b])$.  See \cite{Lee89}.  And, $\mathcal{EBV}([a,b])$
consists of equivalence classes of functions agreeing a.e.  Each
equivalence class of $\mathcal{EBV}([a,b])$ contains a unique function $g$ 
of bounded variation that
is left continuous on $(a,b]$ and right continuous at $a$ such that
$g$ differs from each function in the equivalence class on a set of
measure zero.  If $f$ is in the equivalence class then
$\|f\|_{\mathcal{EBV}}=\inf_h Vh+\|f\|_\infty=Vg+\|g\|_\infty$.  The infimum
is taken over all $h\in\mathcal{BV}([a,b])$ such that $f=h$ a.e.
Note that here $\|f\|_\infty$ is the essential supremum of $f$ and
$\|g\|_\infty$ reduces to the supremum of $g$.
\begin{theorem}\label{T3.4}
(a) $\mathcal A_D([a,b])$ is not complete.
(b) $\mathcal D^{lc}_0([a,b])$ is dense in $D([a,b])$.
(c) The dual space of $\mathcal A_D([a,b])$ is isometrically
isomorphic to $\mathcal{IBV}_b([a,b])$ and to $\mathcal{BV}([a,b])$.
\end{theorem}

{\bf Proof}. (a)
Define
$F_n=\sum_{m=2}^n(-1)^m\chi_{(-1/m,-1/(m+1)]}$.  Then $F_n\in
\mathcal D^{lc}_0([-1,0])$.  Let $F=\sum_{m=2}^\infty(-1)^m
\chi_{(-1/m,-1/(m+1)]}$.  Then $F\in
\mathcal D([-1,0])$ since $F\in L^1([-1,0])$.  Note that
$\|F-F_n\|_A=1/(n+1)-1/(n+2)\to 0$ as $n\to\infty$.  So $F_n\to F$ in 
the Alexiewicz norm.  But $\lim_{x\to 0^-}F(x)$ does not exist so $F$
is not left regulated on $[-1,0]$.  Since $\mathcal D^{lc}_0([a,b])$ and 
$\mathcal A_D([a,b])$ are isomorphic and isometric (Theorem~\ref{T6}) it 
follows that although $\mathcal A_D([a,b])$ is a normed linear space
it is not a Banach space.  

(b)  Let $f\in D([a,b])$ with primitive $F(x)=\int_a^xf(t)\,dt$.
Let $\epsilon>0$.  Find $g\in\mathcal D^{lc}_0([a,b])$ such that
$\|f-g\|_A<\epsilon$. Let $G(x)=\int_a^x g$.  Then 
$\|f-g\|_A =\sup_{I\subset[a,b]}|\int_I(f-g)|\leq 2\|F-G\|_\infty$.
Since $F$ has limits at $a+$ and $b-$ there are
$\alpha, \beta\in\mathbb R$ with $a<\alpha<\beta<b$ such that
if $x\in[a,\alpha]$ then
$|F(x)|<\epsilon/3$ and if $x\in[\beta,b]$ then $|F(x)-F(b)|<\epsilon/3$.
Since $F$ is continuous on $[a,b]$, by
the Weierstrass approximation theorem, there is a polynomial $P$
such that $\|(F-P)\chi_{[\alpha,\beta]}\|_\infty<\epsilon/3$.
Define 
\begin{equation}
G(x)=\left\{\begin{array}{cl}
0, & a\leq x\leq\alpha/2\\
2P(\alpha)(x/\alpha-1/2), & \alpha/2\leq x\leq \alpha\\
P(x), & \alpha\leq x\leq \beta\\
F(\beta), & \beta\leq x\leq b.
\end{array}\right.
\end{equation}
For $\alpha/2\leq x\leq \alpha$ we have
\begin{equation}
|F(x)-G(x)| \leq |F(x)|+ |P(\alpha)|
\leq |F(x)|+ |F(\alpha)-P(\alpha)|+|F(\alpha)|
<\epsilon.
\end{equation}
Hence, $\|F-G\|_\infty<\epsilon$ and $G$ is continuous on $[a,b]$ with a 
continuous derivative
except perhaps at $\alpha/2$, $\alpha$ or $\beta$.  Let
$g(\alpha/2)=0$, $g(\alpha)=G'(\alpha-)$ and $g(\beta)=G'(\beta-)$.
For other values of $x$ let $g(x)=G'(x)$.
Then $g\in\mathcal D^{lc}_0([a,b])$ and
$\|f-g\|_A\leq 2\|F-G\|_\infty<2\epsilon$.  

(c)  The dual spaces of $\mathcal A_D([a,b])$ and $\mathcal D^{lc}_0([a,b])$
are isometrically isomorphic. It is an elementary result of functional analysis that if $Y$ is a
Banach space and $X$ is a dense subspace then $X^\ast=Y^\ast$.  For
example, \cite[p.~194]{KF75}.  By (b) then, the dual space of $\mathcal A_D([a,b])$ 
is 
isometrically isomorphic to the dual of $D([a,b])$.  But this is known
to be $\mathcal{EBV}([a,b])$ \cite{Lee89}.
Now,
the (Lebesgue) integral provides a linear isometry and
isomorphism between $\mathcal{IBV}_b([a,b])$ 
and $\mathcal{EBV}([a,b])$.  If $g\in\mathcal{IBV}_b([a,b])$ and 
$g(x)=\int_a^xh(t)\,dt$ for $h\in\mathcal{EBV}([a,b])$ then
$\|g\|_{\mathcal{IBV}}=\|g'\|_{\mathcal{EBV}}=\|h\|_{\mathcal{EBV}}$.  
The proof of (c) now follows.
\qed

If $T$ is a continuous linear functional on $D([a,b])$ then there is a
unique $g\in\mathcal{EBV}([a,b])$ such that $\langle T,f\rangle=
\int_a^b f(t)g(t)\,dt$
for all $f\in D([a,b])$.  See \cite{Lee89}.
Of course $g$ can be changed on a set of measure zero without affecting
the value of the integral.
The continuous linear functionals on $\mathcal A_D([a,b])$ can now
be characterised.  If $S\in\mathcal A_D^\ast([a,b])$ then there is a
unique function $h\in\mathcal{IBV}_b([a,b])$ such that if $f\in
\mathcal A_D([a,b])$ then $\langle S,f\rangle=
\int_a^bfh=-\int_a^bF(t)h'(t)\,dt$,
where $F\in\mathcal D^{lc}_0([a,b])$ is the primitive of $f$.
Observe that if $(f_n)\subset\mathcal A_D([a,b])$ such that
$\|f_n\|_A\to 0$ then Theorem~\ref{theorem24.2}(a) shows
$\int_a^bf_nh\to 0$, i.e., continuity of the linear functional.

\section {The $LL$ primitive integral and basic properties}\label{S3}
\setcounter{equation}{0}

In this section we will study integrability of distributions whose primitives are left regulated and Lebesgue integrable, and define an integral for such distributions. Properties of the integral, integrable distributions and their primitives are studied.  
  
\subsection{$LLP$ integrability and the $LL$ primitive integral}\label{S31}
We denote
\begin{equation}\label{E130}
\begin{aligned}
\mathcal L^{lr}(I) =&\{H : I\to \mathbb R | \hbox{$H$ is left regulated, locally Lebesgue integrable},\\
&\hbox{and has the right limit at $\inf I$}\},\\
\mathcal L^{lc}(I) =&\{G\in \mathcal L^{lr}(I) | \hbox{$G$ is left continuous, and $G(\min I) = G(\min I+)$ if $\min I$ exists}\},\\
\mathcal L^{lc}_0(I)=&\{F\in \mathcal L^{lc}(I)| F(\inf I+)=0\}.
\end{aligned}
\end{equation}
We use the same notations also in the case when local Lebesgue integrability is replaced by Lebesgue integrability.
When Lebesgue integrability is needed we mention it. 
Properties of the Lebesgue integral ensure that to each  $H \in \mathcal L^{lr}(I)$ there corresponds a unique distribution on $I$, denoted also by  $H$, and defined by (\ref{E0011}),   
where $\smallint$ denotes the Lebesgue integral. 

A distribution $f$ on $I$ is called $LLP$ integrable on $I$ 
if it is the distributional derivative of some primitive $H \in \mathcal L^{lr}(I)$.
 Denote
\begin{equation}\label{E321}
\mathcal A_L(I) = \{f \in \mathcal D' | f \ \hbox{ is $LLP$ integrable on } \ I\}.
\end{equation} 
$\mathcal A_L(I)$ is a linear subspace of $\mathcal A_D(I)$.
If $f$ is $LLP$ integrable with a primitive $H\in \mathcal L^{lr}(I)$, then (\ref{E011})  determines  primitives $G\in \mathcal L^{lc}(I)$ and  $F\in \mathcal L^{lc}_0(I)$ of $f$. 

If $f\in \mathcal A_L(I)$, $F$ is a primitive of $f$ in $\mathcal L^{lc}(I)$  and $a,\,b\in I$, define the $LL$ primitive integral of $f$ from $a$ to $b$ by (\ref{E22}). 
The so obtained integral is unique, additive over intervals, linear, and changes its sign if the integration limits are 
reversed. Proofs of these properties are same as the proofs presented in Theorem \ref{T2} for corresponding properties for the $LD$ primitive integral. 
The fundamental theorem of calculus holds, i.e., the result of Theorem \ref{T3} holds when $f\in \mathcal A_L(I)$ and $G\in \mathcal L^{lc}(I)$, and the proof is same. Thus the mapping $\mathcal F$ defined by
\begin{equation}\label{E223} 
\mathcal F(f)=F, \ f\in\mathcal A_L(I),\ \hbox{where $F$ is the primitive of $f$ in}\ \mathcal L^{lc}_0(I),
\end{equation}
is a linear isomorphism from $\mathcal A_L(I)$ to $\mathcal L^{lc}_0(I)$. The Hake theorem, i.e., Theorem \ref{T4} holds as well for $LL$ primitive integral.

\subsection{Order and norm properties}\label{S33}

Also in $\mathcal L^{lc}_0(I)$, there is the pointwise partial order: $F \le G$ if and only if $F(x) \le G(x)$ for all
$x \in I$. 
This allows us to define a   partial order $\preceq$ in $\mathcal A_L(I)$ by
$f \preceq g$ in $\mathcal A_L(I)$ if and only if  $\mathcal F(f)\le \mathcal F(g)$ in $\mathcal L^{lc}_0(I)$. 
In particular (\ref{E30}) holds, and 
$f\preceq g$ in $\mathcal A_L(I)$ if and only if $F \le G$, where $F$ and $G$ are the
respective primitives in $\mathcal L^{lc}_0(I)$. The proofs that $\mathcal L^{lc}(I)$ and $\mathcal L^{lc}_0(I)$ are ordered vector spaces, and that $\mathcal A_L(I)$ is  order isomorphic to $\mathcal L^{lc}_0(I)$, are the same 
as the proofs presented in Theorem \ref{T5} for corresponding properties for $\mathcal D^{lc}(I)$, $\mathcal D^{lc}_0(I)$ and $\mathcal A_D(I)$.

Next we will show that both $\mathcal L^{lc}_0(I)$ and $\mathcal A_L(I)$ are lattice-ordered vector spaces.
Recall that an ordered vector  space $E$ 
is lattice-ordered if the partial order $\preceq$ of $E$ satisfies the following condition. 
\newline
(3) $x \vee y$ and $x \wedge y$ are in $E$. The join is $x \vee y = \sup\{x, y\} = w$ such that
$x \preceq w$, $y \preceq w$ and if $x \preceq \tilde w$ and $y \preceq \tilde w$ then $w \preceq \tilde w$. The meet is $x \wedge y =
\inf\{x, y\} = w$ such that $w \preceq x$, $w \preceq y$ and if $\tilde w \preceq x$ and $\tilde w \preceq y$ then $\tilde w \preceq w$.

The definitions $(F \vee G)(x) = \sup(F,G)(x) =\max(F(x),G(x))$,
and $(F\wedge G)(x) = \inf(F,G)(x) = \min(F(x),G(x))$ define lattice operations
 in $\mathcal L^{lc}(I)$ and in $\mathcal L^{lc}_0(I)$.

\begin{theorem}\label{T35} (a) $\mathcal L^{lc}_0(I)$ and $\mathcal L^{lc}(I)$ are lattice-ordered. (b)  $\mathcal A_L(I)$ is lattice-ordered.
\end{theorem}

{\bf Proof}. (a) Let $F,\,G \in \mathcal L^{lc}_0(I)$. Define $\Phi = (F \vee G)$ and $\Psi = (F \wedge G)$. We
need to prove $\Phi,\,\Psi \in\mathcal L^{lc}_0(I)$. Let $\inf I<c\le\max I$ and prove $\Phi$ is left continuous at $c$.
Suppose $F(c)>G(c)$. Given $\epsilon > 0$ there is $\delta > 0$ such that $|F(x)-F(c)|<\epsilon$,
$|G(x)-G(c)|<\epsilon$ and $F(x)>G(x)$ whenever $x \in (c-\delta,c)$. For such $x$, $|\Phi(x)- \Phi(c)| = |F(x)-F(c)|<\epsilon$. If $F(c) = G(c)$, then $|\Phi(x)-\Phi(c)| \le \max(|F(x)- F(c)|, |G(x) - G(c)|)<\epsilon$. Therefore, $\Phi$ is left continuous on $I\setminus\{\inf I\}$. 
Similarly, $\Phi$ has the right limit at $\inf I$, and is right continuous at the possible left end point of $I$, so that $\Phi\in\mathcal L^{lc}_0(I)$. Similarly with the
infimum. Hence, $\Phi,\,\Psi\in\mathcal L^{lc}_0(I)$. The proof that $\mathcal L^{lc}(I)$ is lattice-ordered is similar. 

(b) First, we show that $\mathcal A_L(I)$ is closed under the operations $f \vee g$ and
$f \wedge g$. For $f,\, g \in\mathcal A_L(I)$, we have $f \vee g = \sup(f, g)$. There is $h$ such that $f \preceq h$,
$g \preceq h$, and if $f \preceq \tilde h$, $g \preceq \tilde h$, then $h \preceq \tilde h$. This last statement is equivalent to
$\mathcal F(f) \le \mathcal F(g)$, $\mathcal F(g) \le \mathcal F(h)$, and if $\mathcal F(f) \le \mathcal F(\tilde h)$ , $\mathcal F(g) \le \mathcal F(\tilde h)$ , then $\mathcal F(h) \le \mathcal F(\tilde h)$ . But then $\mathcal F(h) = \max(\mathcal F(f),\mathcal F(g))$
and $h = \mathcal F(h)'$ so $f \vee g = (\mathcal F(f) \vee \mathcal F(g))'\in\mathcal A_L(I)$. Similarly, $f \wedge g = (\mathcal F(f) \wedge \mathcal F(g))' \in\mathcal A_L(I)$.
\qed

Let $E$ be a lattice-ordered normed space.
Define  $|x| = x \vee (-x)$, $x^+=x\vee 0$ and $x^-=(-x)\vee 0$. 
Then $x=x^+-x^-$ and $|x|=x^++x^-$.  $E$ is called a normed Riesz space if the norm $\|\cdot\|$ of $E$ is a Riesz-norm, i.e., it satisfies the following condition.
\newline
(4) $|x| \preceq |y|$ implies $\|x\| \le \|y\|$.

Assume next that the functions of $\mathcal L^{lc}(I)$ and $\mathcal L^{lc}_0(I)$ are Lebesgue integrable. We will see that these spaces can be normed by 
 the $L^1$-norm:  
\begin{equation}\label{E310}
\|F\|_1=\int_I|F(t)|\,dt.
\end{equation}
Using the isomorphism $\mathcal F$ we define a 1-norm in $\mathcal A_L(I)$ by $\|f\|_1=\|\mathcal F(f)\|_1$, i.e., %
\begin{equation}\label{E311}
\|f\|_1= \int_I|F(t)|\,dt, \quad  \hbox{ where }\ F=\mathcal F(f).
\end{equation}
The lattice operations show that the $LL$ primitive is absolute: if $f$ is
integrable so is $|f|$.

\begin{theorem}\label{T61} 
Assume that the functions of $\mathcal L^{lc}(I)$ and $\mathcal L^{lc}_0(I)$ are Lebesgue integrable. (a) $\mathcal L^{lc}_0(I)$ and $\mathcal L^{lc}(I)$ are normed Riesz spaces with respect to pointwise ordering and $L^1$-norm. 
(b) $\mathcal A_L(I)$ is a  normed Riesz space with respect to partial order and norm defined by (\ref{E30}) and (\ref{E311}).
(c) If $f\in\mathcal A_L(I)$ then $|f|$, $f^+$ and $f^-$ are in 
 $\mathcal A_L(I)$; $\int_If=\int_If^+-\int_If^-$; $|\int_I f|\geq
|\int_I|f||$.
\end{theorem}

{\bf Proof}. (a) To prove that the $L^1$-norm is a norm in $\mathcal L^{lc}(I)$ and in $\mathcal L^{lc}_0(I)$, 
 first note they are linear subspaces of the space $L^1(I)$ of all Lebesgue integrable functions from $I$ to $\mathbb R$.
And, if $F\in \mathcal L^{lc}(I)$ such that $\|F\|_1 =0$, then 
then $F(x) = 0$ for almost all $x \in I$. Since $F$ is left continuous in $I\setminus\{\inf I\}$, and right continuous at the possible minimum of $I$, then $F(x) = 0$ for all $x \in I$ (see the proof of Theorem \ref{T6}). Positivity, homogeneity
and the triangle inequality are inherited from $L^1(I)$.

The cone $L^1(I)_+=\{G\in L^1(I)| G(x)\ge 0$ for a.e. $x\in I\}$ is closed with respect to the $L^1$-norm. Hence, if $\|F_n-F\|_1\to 0$ in $\mathcal L^{lc}(I)$ and $F_n\in \mathcal L^{lc}(I)_+=\{G\in\mathcal L^{lc}(I)|G(x)\ge 0$ for all $x\in I\}$, then $F(x)\ge 0$ for a.e. $x\in I$. One-sided continuity properties of $F$ ensure that 
$F$ cannot have negative values, whence $F\in \mathcal L^{lc}(I)_+$. This proves that $\mathcal L^{lc}(I)_+$ is closed in
the $L^1$-norm topology of $\mathcal L^{lc}(I)$, so $\mathcal L^{lc}(I)$ is an ordered normed space. The proof that $\mathcal L^{lc}_0(I)$ is an ordered normed space is similar.  

If $|F| \le |G|$, then $\|F\|_1 \le\|G\|_1$. Hence, the $L^1$-norm is a Riesz norm. Thus  $\mathcal L^{lc}(I)$ and 
$\mathcal L^{lc}_0(I)$ are normed Riesz spaces.

(b) Because $\mathcal F$, defined by (\ref{E223}), is an order isomorphism from $\mathcal A_L(I)$ to $\mathcal L^{lc}_0(I)$,  and  $\|f\|_1=\|\mathcal F(f)\|_1$ for all $f\in\mathcal A_L(I)$, then $\mathcal A_L(I)$ is an ordered normed space because $\mathcal L^{lc}_0(I)$ is.
And, if $|f| \preceq | g|$ then $|\mathcal F(f)|' \preceq |\mathcal F(g)|'$ so $|\mathcal F(f)| \le |\mathcal F(g)|$,
that is, $|\mathcal F(f)(x)| \le |\mathcal F(g)(x)|$ for all $x \in I$. Then $\|f\| = \|\mathcal F(f)\|_1 \le \|\mathcal F(g)\|_1 = \|g\|$. Thus the norm defined by (\ref{E311}) is a Riesz norm. This  concludes the proof that 
$\mathcal A_L(I)$ is a normed Riesz space.

(c) Theorem~\ref{T35}(b) establishes the first part.  If 
$F\in \mathcal L^{lc}_0(I)$ is the primitive of $f$ then
$|\int_a^bf|=|F(b)-F(a)|\geq| |F(b)|-|F(a)||=|\int_a^b|F|'| =|\int_a^b|f||$.
\qed

\subsection{Integration by parts, dual space}\label{S34}
The multipliers for
$L^1([a,b])$ are the essentially bounded functions, $L^\infty([a,b])$.
If $f\in L^1([a,b])$ and $g\in L^\infty([a,b])$ then $fg\in L^1([a,b])$.
The dual space of $L^1([a,b])$ is also $L^\infty([a,b])$.  If
$g(x)=\int_a^xh(t)\,dt$ for some $h\in
L^\infty([a,b])\cap L^1([a,b])$ then
$\int_a^bf(t)g(t)\,dt=F(b)g(b)-\int_a^bF(t)h(t)\,dt$, where
$F(x)=\int_a^xf(t)\,dt$.  There are
analogues for the
$LL$ 
primitive integral.

\begin{definition}\label{D34.0}
Let $c\in[a,b]$.  Define
\begin{equation}
\Lambda_c([a,b])=\{g:[a,b]\to\mathbb R\mid
g(x)=\int_c^xh(t)\,dt\quad\text{ for some }h\in
L^\infty([a,b])\cap L^1([a,b])\}.
\end{equation}
\end{definition}
The $L^1$ condition is redundant when $[a,b]$ is compact.
Note that functions in ${\Lambda}_c([a,b])$
are Lipschitz continuous and vanish at $c$.
\begin{definition}\label{D34.1}
Let $f\in\mathcal A_L([a,b])$ with primitive $F\in\mathcal L^{lc}_0([a,b])$.
Let $g\in\Lambda_c([a,b])$ for some $c\in[a,b]$.
Define the integration by parts formula
$\int_a^bfg=F(b)g(b)-\int_a^bF(t)g'(t)\,dt$.
\end{definition}
As in Section~\ref{S24}, density of $C^1$ functions in $L^1$ and 
a sequential approach justifies the definition.
\begin{theorem}\label{theorem34.2}
Let $f\in\mathcal A_L([a,b])$ with primitive $F\in\mathcal L^{lc}_0([a,b])$.
Let $g\in\Lambda_a([a,b])$. (a)  Then
$|\int_a^bfg|\leq|F(b)||g(b)|+\|F\|_1\|g'\|_\infty$.
(b) Suppose $(F_n)\subset C([a,b])\cap C^1((a,b))$ such that $F_n(a)=0$,
$F_n(b)=F(b)$ and 
$\|F'_n-f\|_1\to 0$.
Let $H(x)=
F(x)g(x)-\int_a^xF(t)g'(t)\,dt$ and let 
$H_n(x)=\int_a^xF'_n(t)g(t)\,dt$.
Then
$|H(b)-H_n(b)|
\leq \|F-F_n\|_1\|g'\|_\infty\to 0$ 
as $n\to\infty$.
\end{theorem}

{\bf Proof}. (a) This is the H\"older inequality. 

(b) See the proof of Theorem~\ref{theorem24.2}.
\qed

\begin{proposition}\label{prop34.3}
Let $[a,b]$ be a compact interval.
Suppose $(F_n)\subset C([a,b])\cap C^1([a,b])$ such that
$\|F'_n-f\|_1\to 0$.
Let $H(x)=
F(x)g(x)-\int_a^xF(t)g'(t)\,dt$ and let 
$H_n(x)=\int_a^xF'_n(t)g(t)\,dt$.
Then $H, H_n\in\mathcal L^{lc}_0([a,b])$ and
$\|H-H_n\|_1\to 0$.
\end{proposition}

{\bf Proof}.
The proof is similar to that for Proposition~\ref{theorem24.3}.
Now, using the H\"older inequality and the
Fubini--Tonelli theorem,
\begin{equation}
\begin{aligned}
\|H'-H'_n\|_1&\leq \int_a^b|F(x)-F_n(x)||g(x)|\,dx+\int_a^b
|F(t)-F_n(t)||g'(t)|\int_a^t\,dx\,dt\\
&\leq\|F-F_n\|_1\left(\|g\|_\infty+(b-a)\|g'\|_\infty\right).\qed
\end{aligned}
\end{equation}

Note that if $J$ is a subinterval of $[a,b]$ then $\int_Jf$ is not
defined since $\chi_J$ is not Lipschitz continuous.

Due to the isometric isomorphism between $\mathcal A_L([a,b])$
and $L^1([a,b])$, the dual
space of $\mathcal A_L([a,b])$ is isometrically isomorphic to 
$L^\infty([a,b])$.
\begin{theorem}\label{T4.4}
(a) $\mathcal A_L([a,b])$ is not complete.
(b) $\mathcal L^{lc}_0([a,b])$ is dense in $L^1([a,b])$.
(c) The dual space of $\mathcal A_L([a,b])$ is isometrically
isomorphic to $\Lambda_b([a,b])$ and to $L^\infty([a,b])$.
\end{theorem}

{\bf Proof}. (a)
Use the example in Theorem~\ref{T3.4}(a).  Now,
$\|F-F_n\|_1=1/(n+1)$.

(b) (c) These are essentially the same as in Theorem~\ref{T3.4}.
\qed

If $S\in\mathcal A_L^\ast([a,b])$ then there is a
unique function $g\in\Lambda_b([a,b])$ such that if $f\in
\mathcal A_L([a,b])$ then $\langle S,f\rangle=
\int_a^bfg=F(b)g(b)-\int_a^bF(t)g'(t)\,dt$,
where $F\in\mathcal L^{lc}_0([a,b])$ is the primitive of $f$.

\section {The $LR$ primitive integral and basic properties}\label{S4}
\setcounter{equation}{0}

In this section we will study integrability and an integral of distributions whose primitives are left regulated and Riemann integrable. Properties for the integral, integrable distributions and their primitives are derived.  
  
\subsection{$LRP$ integrability and the $LR$ primitive integral}\label{S41}
We denote
\begin{equation}\label{E140}
\begin{aligned}
\mathcal R^{lr}(I) =&\{H : I\to \mathbb R | \hbox{$H$ is left regulated, locally Riemann integrable,}\\
&\hbox{and has a right limit at $\inf I$}\},\\
\mathcal R^{lc}(I) =&\{G\in \mathcal R^{lr}(I) | \hbox{$G$ is left continuous, and $G(\min I) = G(\min I+)$ if $\min I$ exists}\},\\
\mathcal R^{lc}_0(I)=&\{F\in \mathcal R^{lc}(I)| F(\inf I+)=0\}.
\end{aligned}
\end{equation}
It is well-known that to each  $H \in \mathcal R^{lr}(I)$ there corresponds a unique distribution on $I$, denoted also by  $H$, and defined by (\ref{E0011}),   
where $\smallint$ denotes the Riemann integral. 

A distribution $f$ on $I$ is called $LRP$ integrable on $I$ 
if it is the distributional derivative of some primitive $H \in \mathcal R^{lr}(I)$.
 Denote
\begin{equation}\label{E421}
\mathcal A_R(I) = \{f \in \mathcal D' | f \ \hbox{ is $LRP$ integrable on } \ I\}.
\end{equation} 
$\mathcal A_R(I)$ is a linear subspace of $\mathcal A_L(I)$.
If $f$ is $LRP$ integrable with a primitive $H\in \mathcal R^{lr}(I)$, then (\ref{E011})  determines  primitives $G\in \mathcal R^{lc}(I)$ and  $F\in \mathcal R^{lc}_0(I)$ of $f$. 

If $f\in \mathcal A_R(I)$, $F$ is a primitive of $f$ in $\mathcal R^{lc}(I)$  and $a,\,b\in I$, define the $LL$ primitive integral of $f$ from $a$ to $b$ by (\ref{E22}). 
The so obtained integral is unique, additive over intervals, linear, and changes its sign if the integration limits are 
reversed. Proofs of these properties are the same as the proofs presented in Theorem \ref{T2} for corresponding properties for the $LD$ primitive integral. 
The fundamental theorem of calculus holds, i.e., the result of Theorem \ref{T3} holds when $f\in \mathcal A_R(I)$ and $G\in \mathcal R^{lc}(I)$, and the proof is the same. The mapping
$\mathcal F$, defined by   
\begin{equation}\label{E323} 
 \mathcal F(f)=F, \ f\in\mathcal A_R(I),\ \hbox{where $F$ is the primitive of $f$ in}\ \mathcal R^{lc}_0(I),
\end{equation}
is a linear isomorphism from $\mathcal A_R(I)$ to $\mathcal R^{lc}_0(I)$. The Hake theorem, i.e., Theorem \ref{T4} holds as well for $LL$ primitive integral.

\subsection{Order and norm properties}\label{S430}

Also in $\mathcal R^{lc}_0(I)$, there is the pointwise partial order: $F \le G$ if and only if $F(x) \le G(x)$ for all
$x \in I$. 
We can thus define a   partial order $\preceq$ in $\mathcal A_R(I)$ by
$f \preceq g$ in $\mathcal A_R(I)$ if and only if  $\mathcal F(f)\le \mathcal F(g)$ in $\mathcal R^{lc}_0(I)$. 
In particular (\ref{E30}) holds, and 
$f\preceq g$ in $\mathcal A_R(I)$ if and only if $F \le G$, where $F$ and $G$ are the
respective primitives in $\mathcal R^{lc}_0(I)$. The proofs that $\mathcal R^{lc}(I)$ and $\mathcal R^{lc}_0(I)$ are ordered vector spaces, and that $\mathcal A_R(I)$ is  order isomorphic to $\mathcal R^{lc}_0(I)$, are the same 
as the proofs presented in Theorem \ref{T5} for corresponding properties for $\mathcal D^{lc}(I)$, $\mathcal D^{lc}_0(I)$ and $\mathcal A_D(I)$.

The lattice operations
are defined for $F,\,G \in\mathcal R^{lc}(I)$ by $(F \vee G)(x) = \sup(F,G)(x) =\max(F(x),G(x))$.
And, $(F\wedge G)(x) = \inf(F,G)(x) = \min(F(x),G(x))$. 
The proof that both $\mathcal R^{lc}_0(I)$ and $\mathcal R^{lc}(I)$ are lattice-ordered, and that  $\mathcal A_R(I)$ is lattice-ordered with respect to the partial order $\preceq$ defined above is same as that given in the proof of Theorem
\ref{T35} to $\mathcal L^{lc}_0(I)$, $\mathcal L^{lc}(I)$ and $\mathcal A_L(I)$.

Assume next that $I$ is compact. By Lemma \ref{L000} Riemann integrability of a left regulated function $F:I\to\mathbb R$ is equivalent to the boundedness of $F$.
Thus we can define  the sup-norm norm $\|\cdot\|_\infty$ in   $\mathcal R^{lc}(I)$ and in  $\mathcal R^{lc}_0(I)$:
\begin{equation}\label{E410}
\|F\|_\infty =\sup_{x\in I}|F(x)|.
\end{equation}
Using the isomorphism $\mathcal F$ we define an Alexiewicz norm in $\mathcal A_R(I)$ by $\|f\|=\|\mathcal F(f)\|_\infty$. Equivalently, %
\begin{equation}\label{E411}
\|f\|= \|F\|_\infty= \sup_{x\in I}|F(x)|, \quad  \hbox{ where }\ F(t)=\int_{\min I}^tf.
\end{equation}
We will show that the spaces $\mathcal R^{lc}(I)$ and  $\mathcal R^{lc}_0(I)$ and  $\mathcal A_R(I)$ are Banach lattices,
i.e., complete normed Riesz spaces.

\begin{theorem}\label{T62} Assume that $I$ is compact. (a) $\mathcal R^ {lc}_0(I)$ and $\mathcal R^{lc}(I)$ are Banach lattices with respect to pointwise ordering and sup-norm. 
(b) $\mathcal A_R(I)$ is a  normed Riesz space with respect to partial order and norm defined by (\ref{E30}) and (\ref{E411}).
(c) If $f\in\mathcal A_R(I)$ then $|f|$, $f^+$ and $f^-$ are in 
 $\mathcal A_R(I)$; $\int_If=\int_If^+-\int_If^-$; $|\int_I f|\geq
|\int_I|f||$.

\end{theorem}

{\bf Proof}. (a) Noticing that the spaces $\mathcal R^{lc}(I)$ and in $\mathcal R^{lc}_0(I)$ are also linear subspaces of $L^\infty(I)$, the proof that they are normed Riesz spaces is similar to that given for $\mathcal L^{lc}(I)$ and in $\mathcal L^{lc}_0(I)$ in Theorem \ref{T61} when $L^1(I)$ is replaced by $L^\infty(I)$. 
 To show $\mathcal R^{lc}(I)$  is complete,
suppose $(F_n)$ is a Cauchy sequence in $\mathcal R^{lc}(I)$. Then $(F_n)$ is a Cauchy sequence
in $L^\infty(I)$ so there is $F \in L^\infty(I)$ such that $\|F -F_n\|_\infty \to 0$. To show $F$ is left
continuous in $I\setminus\{\min I\}$, suppose $c \in I\setminus\{\min I\}$. For $x\in (\min I,c)$ and $n \in\mathbb N$,
$$
\begin{aligned}
|F(c)- F(x)| &\le | F(c)- F_n(c)| + |F_n(c) - F_n(x)| + |F_n(x) - F(x)|\\
&\le 2\|F - F_n\|_\infty + |F_n(c) - F_n(x)|.
\end{aligned}
$$
Given $\epsilon > 0$, fix $n$ large enough so that $\|F - F_n\|_\infty < \epsilon/3$. Then let $x\to c-$.
Hence, $F$ is left continuous on $I\setminus\{\min I\}$. 
We can see that $F$ has a right limit at $c=\min I$ by taking $x,\, y >
c$ and letting $x,\, y \to c+$ in $|F(x) - F(y)| \leq 2\|F - F_n\|_\infty + |F_n(x) - F_n(y)|$.
Therefore,
$F \in \mathcal R^{lc}(I)$  and the space $\mathcal R^{lc}(I)$ is complete. The space $\mathcal R^{lc}_0(I)$ is complete since it is a closed subspace of $\mathcal R^{lc}(I)$.  

(b) Because $\mathcal F$ is an order isomorphism from $\mathcal A_R(I)$ to $\mathcal R^{lc}_0(I)$,  and  $\|f\|=\|\mathcal F(f)\|_\infty$ for all $f\in\mathcal A_R(I)$, the proof that $\mathcal A_L(I)$ is a  normed Riesz space is the same as that given for $\mathcal A_L(I)$ in the proof of Theorem \ref{T61}.
To prove it is complete, suppose
$(f_n)$ is a Cauchy sequence in $\mathcal A_R(I)$. Then $\|\mathcal F(f_n) -\mathcal F(f_m)\|_\infty = \|f_n -f_m\|$ so $(\mathcal F(f_n))$
a Cauchy sequence in $\mathcal R^{lc}_0(I)$. There is $F \in \mathcal R^{lc}_0(I)$  such that $\|\mathcal F(f_n) - F\|_\infty \to 0$. And
then $\|f_n - F'\| = \|\mathcal F(f_n) - F\|_\infty \to 0$. Since $F \in\mathcal R^{lc}_0(I)$, we have $F' \in\mathcal A_R(I)$  and 
$\mathcal A_R(I)$ is complete.

(c) The proof is the same as in Theorem~\ref{T61}.
\qed

\subsection{Banach algebra}\label{S44} 
In this subsection we assume that $I$ is compact.

The spaces of Lebesgue and $HK$ integrable functions are not closed under
pointwise multiplication.  For example, if $f(x)=x^{-2/3}$ then
$f$ is Lebesgue integrable on $[0,1]$ but $f^2$ is not.  However,
each of the spaces $\mathcal R^{lr}(I)$, $\mathcal R^{lc}(I)$ and 
$\mathcal R^ {lc}_0(I)$ is closed under pointwise multiplication.
This makes them into commutative Banach algebras.  The  isomorphism between
$\mathcal R^{lc}_0(I)$ and $\mathcal A_R(I)$ makes this latter
space into a commutative Banach algebra.

A commutative algebra is a vector space $V$ over scalar field $\mathbb R$
with a multiplication $V\times V\mapsto V$ such that for all
$u,v,w\in V$ and all $a\in\mathbb R$, $u(vw) =(uv)w$ (associative),
$uv=vu$ (commutative),
$u(v+w)=uv+uw$ and $(u+v)w=uw+vw$ (distributive),
$a(uv)=(au)v$. 
If $(V,\|\cdot\|_V)$ is a Banach space and $\| uv\|_V\leq
\| u\|_V\| v\|_V$ then it is a Banach algebra.

The spaces  $\mathcal R^{lr}(I)$, $\mathcal R^{lc}(I)$ have unit $e=1$.
A unit for $\mathcal R^ {lc}_0(I)$ would need to equal $1$ on 
$(\min I,\max I]$ and $0$ at $\min I$.  But such a function is not right
continuous at $\min I$ so it is not in $\mathcal R^ {lc}_0(I)$.
However, $\mathcal R^ {lc}_0(I)$ has an approximate identity.
If $\min I>-\infty$ define the sequence of continuous functions
$$
u_n(x)=\left\{\begin{array}{cl}
n(x-\min I), & \min I \leq x\leq \min I+1/n\\
1, & \min I +1/n\leq x\leq \max I.
\end{array}
\right.
$$
For each $F\in\mathcal R^{lc}_0(I)$ we have
$\| F-u_nF\|_\infty\to 0$.
$\mathcal R^{lc}_0(I)$ is then said to have an
approximate identity.  A similar construction 
was used in \cite{ET12a} when $I=[-\infty,\infty]$.

\begin{theorem}\label{T63} (a) $\mathcal R^ {lc}_0(I)$ and $\mathcal R^{lc}(I)$ are commutative Banach algebras with respect to pointwise multiplication and sup-norm. 
(b) For $f,g\in\mathcal A_R(I)$, with respective primitives
$F,G\in\mathcal R^{lc}_0(I)$, define a product by $fg=(FG)'$.  Then
$\mathcal A_R(I)$ is a commutative Banach algebra.  It has no unit.
\end{theorem}

{\bf Proof}. (a) Limits and pointwise products commute.  If two functions
vanish at a point so does their product.
 
(b) This follows from the isomorphism between $\mathcal R^ {lc}_0(I)$ and
$\mathcal A_R(I)$.  See \cite[Lemma~1]{ET12a}.
\qed

\subsection{Integration by parts}\label{S45}
The Riemann--Stieltjes integral $\int_a^bF\,dg$ exists for all $F\in
C([a,b])$ and all $g\in\mathcal{BV}([a,b])$. It is defined with a
globally fine partition.  The integral equals $A\in\mathbb R$ if for all
$\epsilon>0$ there is $\delta>0$ such that if $a=x_0<x_1<\ldots<x_n=b$
satisfies $\max_{1\leq i\leq n}|x_i-x_{i-1}|<\delta$ and 
$\xi_i$ is any point in $[x_{i-1},x_i]$ then 
$|\sum_{i=1}^n F(\xi_i)[g(x_i)-g(x_{i-1})]-A|<\epsilon$.
For example, see \cite{ML80}.  For the $LR$ primitive integral we
use a similar type of Stieltjes integral with a locally fine countable
partition of left open intervals.  This yields an integration by
parts formula for which the multipliers are right continuous functions
of bounded variation.  While a necessary condition for existence of 
the Riemann--Stieltjes integral is that at each point of $[a,b]$ one
of $F$ and $g$ is continuous (see \cite[10.6]{H63})
for our version we will require that
$F$ is left continuous and $g$ is right continuous.

The following definitions are necessary.
\begin{definition}\label{D45.1}
A left gauge is a mapping $\gamma$ from $(a,b]$ to the intervals in
$(a,b]$ such that for each $y\in(a,b]$ there is $x\in(a,y)$ such that
$\gamma(y)=(x,y]$.
A $\gamma$-fine left partition is a mutually disjoint collection of
intervals $\mathcal P$ such that $\cup_{I\in\mathcal P}I=(a,b]$; 
if $I\in\mathcal P$ then
$I=(x,y]$ for some $a\leq x<y\leq b$ and $(x,y]\subset\gamma(y)$.
\end{definition}
The definition makes sense for all $-\infty\leq a<b\leq\infty$.  If
$[a,b]$ is a compact interval then the gauge can be constructed from 
a positive function $\delta$ as is done in Henstock--Kurzweil integration.
The  intervals in the gauge are then of the form
$(x-\delta(x),x]$.

A collection of disjoint intervals in $\mathbb R$ is necessarily
countable.  Note that a $\gamma$-fine left partition need not be finite.
For example, if for each $x>0$ we have $\gamma(x)\subset(x/2,x]$
then every $\gamma$-fine left partition of
$[-1,1]$ must be denumerable.  Without loss of generality we will assume
each $\gamma$-fine left partition is denumerable.

\begin{definition}\label{D45.2}
Let $F\in\mathcal R^{lc}([a,b])$ and let $g\in\mathcal{BV}([a,b])$.  The
integral $\int_a^bF\,dg=A\in\mathbb R$ if for each $\epsilon>0$ there is
a left gauge $\gamma$ such that if 
$\mathcal P=\{I_i\}_{i=1}^\infty$ is a $\gamma$-fine
left partition of $[a,b]$ then for each $\xi_i\in I_i$ we have
$|\sum_{i=1}^\infty F(\xi_i)[g(y_i+)-g(x_{i}+)]-A|<\epsilon$.  Here,
$I_i=(x_i,y_i]$ for $i\geq 1$.
\end{definition}
The main properties of this integral are in the following theorem.
\begin{theorem}\label{T45.3}
Let $F, F_i\in\mathcal R^{lc}([a,b])$, let $g,g_i\in\mathcal{BV}([a,b])$
and let $\gamma$ be a left gauge on $[a,b]$. (a) There exists a
$\gamma$-fine left partition. (b) The integral is unique.
(c) Let $a_i,b_i\in\mathbb R$.  Then
$\int_a^b(a_1F_1+a_2F_2)\,d(b_1g_1+b_2g_2)=
a_1b_1\int_a^bF_1\,dg_1+a_1b_2\int_a^bF_1\,dg_2+a_2b_1\int_a^bF_2\,dg_1
+a_2b_2\int_a^bF_2\,dg_2$.
(d) $|\int_a^b F\,dg|\leq \|F\|_\infty Vg$.
(e) Suppose each integral $\int_a^bF_n\,dg$ exists and 
$\|H-F_n\|_\infty\to 0$ for some function $H$.  Then $\int_a^b H\,dg$ 
exists and equals $\lim_{n\to\infty}\int_a^bF_n\,dg$.
(f) Let $F$ be a bounded countably stepped function that is left
continuous on $(a,b]$ and right continuous at $a$.  On $(a,b]$ write
$F=\sum_{i=1}^\infty a_i\chi_{I_i}$ where $\{I_i\}$ is a left continuous
partition of $(a,b]$.  Use the notation of Definition~\ref{D45.2}.  
Then $\int_a^b F\,dg=\sum_{i=1}^\infty a_i
[g(y_i+)-g(x_i+)]$.
(g) The integral exists for each  $F\in\mathcal R^{lc}([a,b])$ and 
$g\in\mathcal{BV}([a,b])$.
(h) If $\int_a^cF\,dg$ and $\int_c^b F\,dg$ exist for some $c\in (a,b)$
then $\int_a^bF\,dg$
exists and equals $\int_a^cF\,dg +\int_c^b F\,dg$.
(i) If $\int_a^b F\,dg$ exists then $\int_x^yF\,dg$ exists for all
$(x,y]\subset(a,b]$.  And, $\int_a^bF\,dg=\int_a^cF\,dg +\int_c^b F\,dg$.
\end{theorem}

{\bf Proof}.
(a) This follows from the construction in Lemma~\ref{L000}.
(b) Given two left gauges $\gamma_1$ and $\gamma_2$, define
$\gamma(x)=\gamma_1(x)\cap\gamma_2(x)$.  Then $\gamma$ is
a left gauge.  The usual uniqueness proof for Henstock--Kurzweil
integrals now shows the integral is unique.  For example, see 
\cite[p.~39]{ML80}.
(c)  This follows from the linearity of the approximating series.
(d) $|\sum_{i=1}^\infty F(\xi_i)[g(y_i)-g(x_i)]|\leq \|F\|_\infty Vg$.
(e) Using the inequality in (d), the proof is essentially the same as for 
uniform convergence of
sequences of Riemann integrable functions.  See \cite[p.~84]{ML80}.
(f) Write $I_i=(b_i,c_i]$.  Define a left gauge $\gamma$ so that 
if $x\in I_i$ then $\gamma(x)\subset I_i$.  If $\mathcal P
=\{(\alpha_i,\beta_i]\}$ is 
a $\gamma$-fine left partition then each $I_i$ is a disjoint union
$I_i=\cup\{J_j\mid J_j\subset I_i\}$.  There is a signed Borel measure
$\mu_g$ such that $\mu_g((x,y])=g(y+)-g(x+)$.  Let $\xi_i\in I_i$.  We then have
\begin{align*}
\sum_{i=1}^\infty F(\xi_i)[g(\beta_i+)-g(\alpha_i+)]
&=\sum_{i=1}^\infty a_i\sum_{J_j\subset I_i}[g(\beta_j+)-g(\alpha_j+)]
=\sum_{i=1}^\infty a_i\sum_{J_j\subset I_i}\mu_g(J_j)\\
&=\sum_{i=1}^\infty a_i\mu_g\left(\cup_{J_j\subset I_i}J_j\right)
=\sum_{i=1}^\infty a_i\mu_g\left(I_i\right)\\
&=\sum_{i=1}^\infty a_i[g(c_i+)-g(b_i+)].
\end{align*}
Due to the estimate in (d) all of the series converge absolutely.
(g) By Lemma~\ref{L000}(b) $F$ is the uniform limit of a sequence of
bounded countably stepped functions.  The result now follows from (e)
and (f).
(h) Given $\epsilon>0$ there exist gauges $\gamma_1$ and $\gamma_2$
so that sums over respective $\gamma_i$-fine partitions approximate
$\int_a^cF\,dg$ and $\int_c^bF\,dg$ with error at most $\epsilon/2$.
Defining $\gamma(x)=\gamma_1(x)$ if $x\in(a,c]$ and $\gamma_2(x)$ if
$x\in(c,b]$ defines a left gauge on $(a,b]$ with the property that
sums over $\gamma$-fine partitions approximate
$\int_a^bF\,dg$ with error at most $\epsilon$.
(i) For each $(x,y]\subset(a,b]$, the function $F\chi_{(x,y]}$ is in 
$ R^{lc}([a,b])$.  By (h), $\int_a^b F\chi_{(x,y]}\,dg$ exists.  
Examining the approximating sums shows $\int_a^b F\chi_{(x,y]}\,dg=
\int_x^y F\,dg$.  Now write $\int_a^b F\,dg=\int_a^bF(\chi_{(a,c]}+
\chi_{(c,b]})\,dg$. By linearity this equals $\int_a^c F\,dg+
\int_c^bF\,dg$.
\qed

\section{Examples and Remarks}\label{S6}

In this section we will first construct examples of left regulated functions that are locally integrable in the $HK$, Lebesgue or Riemann sense, and  have at every rational number a discontinuity of the second kind. These functions are also
used to define  primitives which are locally integrable in $I=[0,\infty)$.
Lemmas \ref{L002} and \ref{L003} are used to verify the local integrability of
the constructed functions.
 
\begin{example}\label{Ex601}  Define a mapping $G:\mathbb R\to \mathbb R$ by  
\begin{equation}\label{E47} 
G(t)=\sum_{n=1}^\infty\frac 1{n^2}\left(2(nt-[nt])\cos\left(\frac {\pi}{2(nt-[nt])}\right)+\frac \pi 2\sin\left(\frac {\pi}{2(nt-[nt])}\right)\right), \quad t\in\mathbb R, 
\end{equation}
where $[nt]=m,\ m\le nt< m+1$. For each fixed $m\in\mathbb N$, denote by $G_m(t)$ the $m$th partial sum  of the series (\ref{E47}) when $t\in\mathbb R$. It is easy to verify that the so obtained functions $G_m:\mathbb R\to\mathbb R$ are left regulated, and that the set of all discontinuity points of $G_m$ is 
$$
Z_m=\{\frac ij|\, j\in\{1,\dots,m\}, \ i\in\mathbb Z,\ \hbox{ and $i$ and $j$ are coprime}\}.
$$ 
Moreover, the sequence $(G_m)$ converges uniformly to $G$ on each compact subinterval of $\mathbb R$.  
Define a function $F:\mathbb R\to\mathbb R$ by
\begin{equation}\label{E48}
F(t)=\sum_{n=1}^\infty\frac{(nt-[nt])^2}{n^3}\cos\left(\frac {\pi}{2(nt-[nt])}\right), \quad t\in\mathbb R. 
\end{equation}
The $m$th partial sums of the series (\ref{E48}) define functions $F_m:\mathbb R\to\mathbb R$. Obviously, each $F_m$ is continuous, and the sequence $(F_m)$ converges uniformly to $F$ on each compact subinterval of $\mathbb R$, whence $F$ is continuous. Moreover, $F_m'(t)=G_m(t)$ for each $t\in\mathbb R\setminus Z_m$. Consequently, the hypotheses of 
Lemma \ref{L003} are valid for $F$ and $G$, so that  
 $F'(t)=G(t)$ for each $t\in\mathbb R\setminus\underset{m}{\cup}Z_m$. Thus  
$G$ is by Lemma \ref{L002} locally $HK$ integrable. Because $G$ is locally bounded, it is also locally Riemann integrable. 
$G$ is discontinuous at every point of the set $\underset{m}{\cup}Z_m$, which is the set $\mathbb Q$ of all rational numbers. Moreover, all the discontinuities are of the second kind because of the sine term in the right hand side of (\ref{E47}). 
On the other hand, for each $t\in \mathbb R\setminus\mathbb Q$, the functions $G_m$ are continuous at $t$ and converge uniformly in $[t-1,t+1]$ to $G$, whence $G$ is continuous at $t$. Similarly, since for  every $m\in\mathbb N$, $G_m$ has a left limit at each point of $\mathbb R$, this property holds also for $G$, i.e., $G$ is left regulated.    

The above reasoning shows that (\ref{E47}) defines a  function $G:\mathbb R\to\mathbb R$ that has the following properties:
\begin{itemize}
\item $G$ is left regulated and locally Riemann integrable;
\item $G$ is continuous
in $\mathbb R\setminus Q$, and  each point of $\mathbb Q$ is its discontinuity point of the second kind.
\end{itemize} 

The function $t\mapsto tG(t)$ has the above properties, and it is right continuous at the origin. Its restriction to $I=\mathbb R_+$ belongs to $\mathcal R^{lr}(I)$. The function $G_0(t)=\begin{cases} tG(t), \ t\in I\setminus\mathbb Q_+,\\tG(t-), \ t\in\mathbb Q_+,\end{cases}$  belongs to $\mathcal R_0^{lc}(I)$.   

Also the function $t\mapsto e^{-|t|}G(t)$ has  the properties listed above, and it belongs to $\mathcal R^{lr}(\mathbb R)$. Moreover, it is $HK$ integrable. 
\end{example}

In the next example we present locally Lebesgue integrable  primitives  that are discontinuous at every rational point of their domains, 
and are not locally Riemann integrable.

\begin{example}\label{Ex602} Let $G$ and $F$ be defined by (\ref{E47}) and (\ref{E48}). Define functions  $G^m:\mathbb R\to \mathbb R$, $m\in\mathbb N$, by  
\begin{equation}\label{E409}
G^m(t)=G(t)+\sum_{n=1}^m\frac{1}{2\sqrt{nt-[nt]}}, \quad t\in\mathbb R. 
\end{equation}
It is elementary to verify that $\mathbb Q$ is the set of discontinuity points of functions $G^m$, and that these functions are left regulated. Define functions $F^m:\mathbb R\to\mathbb R$, $m\in\mathbb N$, by
\begin{equation}\label{E600}
F^m(t)=F(t)+\sum_{n=1}^m\frac{[nt]+\sqrt{nt-[nt]}}{n}, \quad t\in\mathbb R_-. 
\end{equation}
$F^m$ is continuous, and $(F^m)'(t)=G^m(t)$ for all $t\in\mathbb R\setminus\mathbb Q$ and $m\in\mathbb N$. Lemma \ref{L002} implies then that functions $G^m$ are locally $HK$ integrable. Because the functions $F^m$ are locally absolutely continuous, then every $G^m$ is locally Lebesgue integrable. But $G^m$ is not locally bounded, and hence not locally Riemann integrable, for any $m\in\mathbb N$.

The functionm $t\mapsto tG^m(t)$ have the above properties, and they are right continuous at the origin. Their restrictions to $\mathbb R_+$ belong to $\mathcal L^{lr}(\mathbb R_+)$, but not to $\mathcal B^{lr}(\mathbb R_+)$. The functions $t\mapsto tG^m(t-)$, restricted to $\mathbb R_+$,  belong to $\mathcal L_0^{lc}(\mathbb R_+)$, but not to $\mathcal R^{lc}_0(\mathbb R_+)$.   

The functions $t\mapsto e^{-|t|}G^m(t)$ are Lebesgue integrable on $\mathbb R$.
\end{example}

Locally $HK$ integrable primitives which are discontinuous at every rational point of their domains, and are not locally Lebesgue integrable, are presented in the next example.

\begin{example}\label{Ex603} Let $G$ and $F$ be defined by (\ref{E47}) and (\ref{E48}). Define functions $G_m:\mathbb R\to \mathbb R$, $m\in\mathbb N$, by  
\begin{equation}\label{E407}
G_m(t)=G(t)+\sum_{n=1}^m\left(\cos\left(\frac{\pi}{2(nt-[nt])}\right) +\frac{\pi\sin(\frac{\pi}{2(nt-[ nt])})}{2(nt-[nt])}\right), \quad t\in\mathbb R. 
\end{equation}
$G_m$ is left regulated, and $\mathbb Q$ is the set of its discontinuity points. Functions $F_m:\mathbb R\to\mathbb R$, defined by
\begin{equation}\label{E408}
F_m(t)=F(t)+\sum_{n=1}^m\frac{1}{n}(nt-[nt])\cos\left(\frac{\pi}{2(nt-[nt])}\right), \quad t\in\mathbb R, 
\end{equation}
are continuous, and $F_m'(t)=G_m(t)$ for all $t\in\mathbb R\setminus\mathbb Q$. 
It then follows from Lemma \ref{L002} that the functions $G_m$ are locally $HK$ integrable. On the other hand, $G_m$ is neither locally Lebesgue integrable nor locally Riemann integrable for any $m\in\mathbb N$, since $F_m$ is not locally absolutely continuous, and $G_m$ is not locally bounded for any $m\in\mathbb N$.

The functions $t\mapsto t^2G_m(t)$, $t\in \mathbb R+$, belong to $\mathcal D^{lr}([0,1])\setminus \mathcal L^{lr}([0,1])$, and the functions  $t\mapsto t^2G_m(t-)$, $t\in\mathbb R_+$,  belong to $\mathcal D^{lc}_0([0,1])\setminus \mathcal L^{lc}_0([0,1])$.

The functions $t\mapsto e^tG_m(t)$ are $HK$ integrable on $\mathbb R$.
\end{example}

An example of a left regulated function that is not  $HK$ integrable at any subinterval of $\mathbb R$  that contains origin is
$$
G_p(t)=\begin{cases} G(t)+\frac{1}{t},\ &t> 0,\\ G(t), \ &t\le 0,\end{cases}
$$ 
where $G$ is defined by (\ref{E47}).

$LCP$ integrable distributions are defined in \cite{SH122}. The space of their primitives is
$$
\mathcal B^{lr}(I)=\{H:I\to \mathbb R|H \hbox{ is bounded and left regulated, and $H(\min I+)$ exists}\}.
$$
$LCP$ integrable distributions have unique primitives in the space
$$
\mathcal B_0(I)=\{F:I\to \mathbb R|F \hbox{ is bounded and left continuous, and $F(\min I+)=F(\min I)=0$}\}.
$$
Because a left regulated function is by Lemma \ref{L001} locally Riemann integrable if and only if it is locally bounded, then $\mathcal B^{lr}(I)= \mathcal R^{lr}(I)$ and $\mathcal B_0(I)=\mathcal R^{lc}_0(I)$.  

In \cite{SH12}  $RP$-integrable distributions are defined so that their primitives form  the space of  
all regulated functions on $I$.  These distributions have unique primitives in the space 
$$
\mathcal P_R(I)=\{F:I\to\mathbb R|F \ \hbox{is regulated and left continuous, and } \ F(\min I)=F(\min I+)=0\}. 
$$
$P_R(I)$ is a proper subset of $B_0(I)$. 
For instance, the function $G_0$ defined in Example \ref{Ex601}, belongs to $\mathcal B^{lr}(\mathbb R_+)$, but  not to $\mathcal P_R(\mathbb R_+)$.
The restriction of the function $F$, defined by (\ref{E611}),  to any interval $I$ for which $\min I$ exists and is irrational belongs to  $\mathcal P_R(I)$. 

To construct  regulated functions which are discontinuous at each rational point, let
$p > 1$ be fixed.
Define a function $F:\mathbb R\to\mathbb R$ by
\begin{equation}\label{E611}
F(t)=  \sum_{n=1}^\infty \frac{1+nt-\left\lfloor nt \right\rfloor}{n^p},  \quad t\in\mathbb R,
\end{equation}
where $\left\lfloor nt \right\rfloor=m,\ m-1<nt\le m, \ m=0,1,\dots$. The reasoning used in Example \ref{Ex601} shows that $F$ is well-defined, that the set of discontinuity points of $F$ is formed by all rational numbers, and that $F$ is left continuous and has right limit at each $t\in\mathbb Q$. Thus $F$ is regulated.
The restriction of $F$ to any interval $I$ for which $\min I$ exists and is irrational belongs to $\mathcal P_R(I)$.

In  \cite{ET09} a theory is presented for the regulated primitive integral of distributions whose primitives belong to the space 
$$
\mathcal B_R =\{F :\mathbb R\cup\{\pm\infty\}\to \mathbb R | F \ \hbox{is regulated and left continuous on} \ \mathbb R,\,F(-\infty) = 0,F(\infty)\in\mathbb  R\}.
$$
The function $F_0$, defined by 
$F_0(t)=e^{-|t|}F(t)$, $t\in\mathbb R$, $F_0(\pm\infty)=0$, belongs to $\mathcal B_R$. 

\begin{remark}\label{R6} 
The proofs of Theorems \ref{T2}, \ref{T5}, \ref{T6}, \ref{T35}, \ref{T61} and \ref{T62} are similar to proofs of corresponding results of \cite{ET09}.
The case of regulated primitives considered in this paper yields 
four different integrals over each of the intervals $(a,b)$, $[a,b]$,
$[a,b)$ and $(a,b]$.  For example, $\int_{(a,b)}F'=F(b-)-F(a+)$.  Since
we have left continuity of primitives with the $LD$, $LL$ and $LR$ primitive 
integrals, the 
integral \eqref{E22} can be considered as
$\int_a^bF'=\int_{[a,b)}F'=F(b-)-F(a-)=F(b)-F(a)$.

The above examples  show that inclusions in $\mathcal P_R(I)\subset\mathcal R^{lc}_0(I)\subset\mathcal L^{lc}_0(I)\subset \mathcal D^{lc}_0(I)$ are proper.
This result and the bijective correspondence $f\stackrel{\mathcal F}{\leftrightarrow} F$ between integrable distributions and their primitives implies  
that the $RP$ integrable distributions form a proper subset of $LCP$ integrable distributions, which is equal to the space of $LRP$ integrable distributions. They in turn form a proper subset of $LLP$ integrable distributions, which form 
a proper subset of $LDP$ integrable distributions.

In defining $\mathcal D^{lr}(I)$, $\mathcal L^{lr}(I)$ and $\mathcal R^{lr}(I)$ we have chosen the primitives to be left regulated. Another obvious choice is to take primitives that are right regulated.  One can also use a convex combination of left and right limits.  As for properties of right regulated functions see \cite[Subsection 1.4.2]{HL94}.

The order $\preceq$ defined by (\ref{E30}) is not compatible with the usual
order on distributions: if $T,U \in\mathcal D'$ then $T \ge U$ if and only if $\left\langle T - U,\phi\right\rangle \ge 0$ for all $\phi\in\mathcal D$ such that $\phi\ge 0$. Nor is it compatible with pointwise ordering in
the case of functions in $\mathcal A_D(I)$. For example, if $f(t) = H_1(t) \sin(t^2)$, where $H_1$ is the Heaviside step function,   i.e., $H_1(t)=\begin{cases} 1, t>0,\\ 0, t\le 0,\end{cases}$  then $0\le F$
so $0\preceq f$ in $\mathcal A_D[-1,1])$ but not pointwise. And, $f$ is not positive in the distributional
sense. Note, however, that if $f \in\mathcal A_D(I)$ is a nonnegative function
or distribution then $0\preceq f$ in $\mathcal A_D(I)$.

In Lebesgue and Henstock--Kurzweil integration, we have equivalence
classes of functions that agree almost everywhere. In $\mathcal A_D(I)$, $\mathcal A_L(I)$ and $\mathcal A_R(I)$ there are no such equivalence classes, for two distributions are equal if they agree on all
test functions. 

Except the last remark  we assume from now on that $I$ is a compact real interval. If function $f$ is continuous, or in
$L^p(I)$ for some $1\leq p\leq\infty$, or in the Denjoy space $D(I)$ of $HK$ integrable functions from $I$ to $\mathbb R$, then 
$f\in\mathcal A_R(I)$ and hence in  $\mathcal A_L(I)$ and
$\mathcal A_D(I)$.  If primitive function $F$ is continuous then $F'\in
\mathcal A_R(I)$.  Taking $F$ to be continuous but with a pointwise 
derivative nowhere, we see that $\mathcal A_R(I)$ contains distributions
that have no pointwise values.  Note that the formula $\int_a^bF'=F(b)-F(a)$
holds although this now has no meaning as a Riemann, Lebesgue or 
Henstock--Kurzweil integral.  If $F$ is continuous but $F'(x)=0$
almost everywhere then the Lebesgue integral $\int_a^b F'(x)\,dx=0$ but
in $\mathcal A_R(I)$ we have  $\int_a^bF'=F(b)-F(a)$. 

The space $\mathcal A_R(I)$ contains all signed Borel measures on 
$(\min I,\max I)$.
Suppose $\mu$ is a signed Borel measure such that 
$\mu(\{\min I\})=\mu(\{\max I\})=0$.  Define $F(\min I)=0$ and
$F(x)=\int_{(0,x)}d\mu$ for $x\in(\min I,\max I]$.  This
Lebesgue integral defines primitive $F$ that is of bounded variation 
and in $\mathcal R^{lc}_0(I)$.  Hence, $\mu\in\mathcal A_R(I)$.  For
example, the Dirac measure, $\delta$, is in $\mathcal A_R([a,b])$ for
any $a<0<b$.

Lemma~\ref{L000} shows left regulated functions can be approximated
uniformly by countably stepped functions.  However, countably stepped functions
are not dense in $\mathcal R^{lc}(I)$.  For example, let
$F=\sum_{n=1}^\infty (-1)^n\chi_{((n+1)^{-1},n^{-1}]}$.  Then
$F\in\mathcal R^{lc}([-1,1])$ and $\|F\|_\infty=1$.  Suppose $\sigma$
is a step function.  Then $\lim_{x\to 0^+}\sigma(x)$ exists.  We have
\begin{equation}
\|F-\sigma\|_\infty\geq \limsup_{x\to 0^+}|F(x)-\sigma(x)|=
\limsup_{x\to 0^+}|F(x)-\sigma(0+)|\geq 1.
\end{equation}

The following example shows $\mathcal D^{lc}_0(I)$ is not closed.  Define
$F_n=\sum_{m=1}^n(-1)^m\chi_{(1/(m+1),1/m]}$.  Then $F_n\in
\mathcal D^{lc}_0([0,1])$.  Let $F=\sum_{m=1}^\infty(-1)^m
\chi_{(1/(m+1),1/m]}$.  Then $F\in
\mathcal D([0,1])$ since $F\in L^1([0,1])$.  Note that
$\|F-F_n\|_A=1/(n+1)-1/(n+2)\to 0$ as $n\to\infty$.  So $F_n\to F$ in 
the Alexiewicz norm.  But $\lim_{x\to 0^+}F(x)$ does not exist so $F$
is not regulated on $[0,1]$.  Hence, although $\mathcal D^{lc}_0(I)$ is a normed
space it is not a Banach space.  The same functions show $\mathcal L^{lc}_0(I)$
is not closed.  Now, $$
\|F_n-F\|_1=\sum_{m=n+1}^\infty\left(\frac{1}{m}-\frac{1}{m+1}\right)=
\frac{1}{n+1}\to 0\quad\text{ as } n\to\infty.
$$
Hence, $\mathcal L^{lc}_0(I)$ is also not a Banach space.  
It now follows that neither $\mathcal A_D(I)$ nor $\mathcal A_L(I)$ are
Banach spaces.  Since the continuous functions are dense in $\mathcal D^{lc}(I)$
the completion of $\mathcal D^{lc}(I)$ in the Alexiewicz norm is $D(I)$ and
the completion of $\mathcal A_D(I)$ is the space of distributional derivatives
of $HK$ integrable functions.  This space was studied in \cite{ET08}.  
Similarly, the completion of $\mathcal A_L(I)$
in the $1$-norm is the space of distributional derivatives of Lebesgue
integrable functions.  This space was studied in \cite{ET12}.

Transfinite series whose terms are indexed with inversely well-ordered sets of real numbers are applied  in \cite{SH13}
to derive further integrability criteria for left regulated functions. 
Such series are used also to prove the converse of Lemma \ref{L002} for left regulated functions. Lemma \ref{L002}, combined with that converse, imply a Fundamental Theorem of Calculus for left regulated functions. It is important in applications to ordinary and impulsive differential equations.

\end{remark}

\section{Applications to distributional Cauchy systems}\label{S8} 
\setcounter{equation}{0}

In this section we will study the following system of distributional Cauchy problems:
\begin{equation}\label{E1}    
 y_i'=f_i(y_1,\dots,y_m), \quad  y_i(a)=c_i, \quad i=1,\dots,m. 
 \end{equation}
Dependence of solutions on $f_i$ and $c_i$ is also studied.  Values of the functions $f_i$ are distributions on a half-open real interval $I=[a,b)$, $a < b\le\infty$. 

The regulated primitive integral is studied in detail in \cite{ET09} when $I=\overline{\mathbb R}$, and applied in \cite{SH12} to problem (\ref{E1}) when $I=[a,b]$.  The left continuous primitive integral is applied in \cite{SH122} to  problem (\ref{E1}) when $I=[a,b]$.
Because $\mathcal B_0(I)=\mathcal R^{lc}_0(I)$ when $I$ is compact, the left continuous primitive integral and the $LR$ primitive integral are equal. Therefore we study only applications of the $LD$ primitive integral and the $LL$ primitive integral to problem (\ref{E1}).  No continuity hypotheses are imposed on functions $f_i$.

\subsection{On the smallest and greatest solutions}\label{S81}

We will first study the existence of the smallest and greatest solutions of problem (\ref{E1}) and their dependence on
$f_i$ and $c_i$.
Component functions $y_i:[a,b)\to\mathbb R$ of solutions of (\ref{E1}) are assumed to be in the space
$\mathcal D^{lc}([a,b))$ of those functions from $[a,b)$ to $\mathbb R$ which are locally $HK$ integrable, right continuous
at $a$ and left continuous on $(a,b)$. 

Assume that   $\mathcal D^{lc}([a,b))$ is ordered pointwise, that the space  $\mathcal A_D([a,b))$ is equipped with the ordering $\preceq$ defined by (\ref{E30}), that  the space $HK_{loc}([a,b))$ of locally $HK$ integrable functions on $[a,b)$ is equipped with  a.e. pointwise ordering, and that a.e. equal functions of $HK_{loc}([a,b))$ are identified. The product spaces $\mathcal D^{lc}([a,b))^m=\times_{i=1}^m\mathcal D^{lc}([a,b))$ and $HK_{loc}([a,b))^m=\times_{i=1}^mHK_{loc}([a,b))$
are ordered by componentwise ordering, i.e., if $x=(x_i,\dots,x_m)$ and $y=(y_1,\dots,y_m)$ belong to one of these product spaces, then 
$$
x\le y \ \hbox{ iff $x_i\le y_i$ for all } \ i=1,\dots,m.
$$

\begin{definition}\label{D21} A function
$(y_1,\dots,y_m)\in\mathcal D^{lc}([a,b))^m$ is called a subsolution of (\ref{E1}) if  
\begin{equation}\label{E101}
y_i'\preceq f_i(y_1,\dots,y_m) \ \hbox{
in $\mathcal A_D([a,b))$, and  $y_i(a)\le c_i$ for every} \ i=1,\dots,m.
\end{equation}
If reversed inequalities hold in (\ref{E101}), we say that $(y_1,\dots,y_m)$ is a supersolution of (\ref{E1}).
If equalities hold in (\ref{E101}), then $(y_1,\dots,y_m)$ is called a solution of (\ref{E1}). 
\end{definition}
The following result that transforms the system (\ref{E1}) into a system of integral equations is a direct consequence of Theorem \ref{T3}.

\begin{lemma}\label{L20} Assume that $(y_1,\dots,y_m)\in \mathcal D^{lc}([a,b))^m$, and that $f_i(y_1,\dots,y_m)\in\mathcal A_D([a,b))$ for every $i=1,\dots,m$. Then
$(y_1,\dots,y_m)$ is a solution of the system  (\ref{E1})  
if and only if it is a solution 
of  the following system of integral equations: 
\begin{equation}\label{E020}
    y_i(t)=c_i+\int_a^tf_i(y_1,\dots,y_m), \quad t\in[a,b),\ i=1,\dots,m.
\end{equation} 
\end{lemma}

The application of monotone methods to find solutions of (\ref{E1}) is complicated by the fact that the limit function, supremum and/or infimum of a pointwise
convergent monotone sequence of $\mathcal D^{lc}(I)$ are not necessarily in $\mathcal D^{lc}(I)$ even in the case when the interval $I$ is compact. For instance, the sequence of functions $x_n\in \mathcal D^{lc}([0,1])$, $n=0,1,\dots$, defined by
$$
x_n(t)=\begin{cases} 0, \ & 0\le t \le \frac 1{2^{n+1}},\\ 1-(-1)^k, \ & \frac 1{2^{k+1}}< t\le \frac 1{2^k}, \ k=0,\dots,n,\end{cases} 
$$
is increasing in the pointwise ordering of $\mathcal D^{lc}[0,1])$, but neither its pointwise limit nor its supremum is in $\mathcal D^{lc}[0,1])$. 

Therefore we study in this section the existence of such solutions of the system (\ref{E1}) whose components
are locally $HK$ integrable on $[a,b)$.

\begin{definition}\label{D201} Given partially ordered sets $X=(X,\le)$ and $Y=(Y,\preceq)$, we say that a mapping $f:X\to Y$ is
{\em increasing} if  $f(x)\preceq f(y)$ whenever $x\le y$ in $X$, and {\em order-bounded}
if there exist  $\underline y,\overline y\in Y$ such that the range $f[X]$ of $f$ is contained in the order interval
$[\underline y,\overline y]=\{y\in Y:\underline y\preceq y \preceq \overline y\}$ of $Y$.
\end{definition}

The first existence and comparison theorem for the smallest and greatest solutions of the system (\ref{E1}) reads as follows.  

\begin{theorem}\label{T21} Assume that $f_i:HK_{loc}([a,b))^m\to \mathcal A_D([a,b))$ is increasing, that the system (\ref{E1}) has in $\mathcal D^{lc}([a,b))^m$ a subsolution $\underline y=(\underline y_1,\dots,\underline y_m)$ and a supersolution $\overline y=(\overline y_1,\dots,\overline y_m)$, and that $\underline y_i\le \overline y_i$ for each $i=1,\dots,m$.
Then the system (\ref{E1}) has in the order interval $[\underline y,\overline y]$ of $\mathcal D^{lc}([a,b))^m$ the smallest and greatest solutions, and they are increasing with respect to $f_i$ and $c_i$, $i=1,\dots,m$.
\end{theorem}

{\bf Proof}. Let $F:HK_{loc}([a,b))^m\to \mathcal D^{lc}([a,b))^m$ be defined by 
\begin{equation}\label{E02}
\begin{cases}
F(x)=(F_1(x_1,\dots,x_m),\dots,F_m(x_1,\dots,x_m)), \quad \hbox{where}\\
F_i(x_1,\dots,x_m)(t)=c_i+\int_a^tf_i(x_1,\dots,x_m),  \quad t\in[a,b),\ i=1,\dots,m.
\end{cases}
\end{equation}
In view of Lemma \ref{L20}, $y=(y_1,\dots,y_m)$  is a solution of the system (\ref{E1}) in  $\mathcal D^{lc}([a,b))^m$  if and only if $y$ is a solution of the  fixed point equation 
    $y=F(y)$
in  $\mathcal D^{lc}([a,b))^m$. Assume that $x=(x_1,\dots,x_m)$ belongs to the order interval $[\underline y,\overline y]$ of $HK_{loc}([a,b))^m$. The given hypotheses imply by Definitions \ref{D21} and \ref{D201} that for every $i=1,\dots,m$,
$$
\underline y_i'\preceq f_i(\underline y)\preceq f_i(x_1,\dots,x_m)\preceq f_i(\overline y)\preceq\overline y_i'
$$ 
for all $i=1,\dots,m$. Moreover, $\underline y_i(a)\le c_i\le\overline y_i(a)$ for every $i=1,\dots,m$.
Thus
$$
\begin{cases}
 \underline y_i(t)\le\underline y_i(a)+\int_a^t f_i(\underline y)\le c_i+\int_a^tf_i(x)\le \overline y_i(a)+\int_a^tf_i(\overline y)\le\overline y_i(t), \\ t\in[a,b),\ i=1,\dots,m.
 \end{cases}
 $$
Because $t\mapsto \int_a^tf_i(x)$ belongs to $\mathcal D^{lc}([a,b))$, it is $HK$ locally integrable. Thus 
$F_i(x)$ is Lebesgue measurable and order bounded by functions $\underline y_i$ and $\overline y_i$ of $HK_{loc}([a,b))$. It then follows from \cite[Proposition 9.39 and Remark 9.25]{CH11} that $F_i(x)\in HK_{loc}([a,b))$. This holds for every 
$i=1,\dots,m$, whence $F(x)\in HK_{loc}([a,b))^m$.

The above results imply that
$F$ maps order interval $[\underline y,\overline y]$ of $HK_{loc}([a,b))^m$ into itself. 
Moreover, $F$ is increasing because the functions $f_i$ are increasing.
By \cite[Theorem 1.1.1]{HL94} there is a unique
chain $C$ in $HK_{loc}([a,b))^m$ that is well-ordered (every non-empty subset of $C$ has the smallest element), and that satisfies
\begin{description}
\item[(I)] {$\underline y = \min
C$, \ and if $\underline y < x$, then $x\in C$ \ iff \ $x = \sup F[\{y\in C:y < x\}]$}.
\end{description}
Since $C$ is well-ordered and $F$ is increasing, then  $W = F[C]$ is well-ordered.
 For every $i=1,\dots,m$, the set
$
W_i=\{z_i:(z_1,\dots,z_m)\in W\}
$
is a well-ordered chain in the order interval $[\underline y_i,\overline y_i]$ of $HK_{loc}([a,b))$.
Thus  $y_{*i}=\sup W_i$ exists in $HK_{loc}([a,b))$ by \cite[Proposition 9.39]{CH11}. Obviously,
$(y_{*1},\dots,y_{*m})$ is the supremum of $W=F[C]$ in $HK_{loc}([a,b))^m$. 
It then follows  from  \cite[Theorem 1.2.1]{HL94}  that $y_*=(y_{*1},\dots,y_{*m})=\max C$, and that $y_*$ is the smallest fixed point of $F$ in the order interval $[\underline y,\overline y]$ of  $HK_{loc}([a,b))^m$. Moreover, every fixed point of $F$ belongs to $\mathcal D^{lc}([a,b))^m$. Thus $(y_{*1},\dots,y_{*m})$
 is  the smallest solution of the system (\ref{E1}) in order interval $[\underline y,\overline y]$ of  $\mathcal D^{lc}([a,b))^m$.
    
 According to \cite[Proposition 1.2.1]{HL94} there exists a unique chain $D$ that is inversely well-ordered,  and that satisfies
\begin{description}
\item[(II)] {$\overline y = \max D$, and if $x < \overline y$, then $x\in D$  iff  $x = \inf F[\{y\in D: x< y\}]$.}
\end{description}
 The proof that  $y^* = \inf F[D]=\min D$ exists and is the greatest fixed point of $F$ in the order interval $[\underline y,\overline y]$ of  $\mathcal D^{lc}([a,b))^m$
 is similar to the above proof. Thus $y^*=(y_1^*,\dots,y_m^*)$ is  the greatest solution of  the system (\ref{E1}) in the order interval $[\underline y,\overline y]$ of  $\mathcal D^{lc}([a,b))^m$.
 Moreover, according to  \cite[Theorem 1.2.1 and Proposition 1.2.1]{HL94},
\begin{equation}
y_*=\min\{x\in[\underline y,\overline y]:F(x)\le x\}, \ y^*=\max\{x\in[\underline y,\overline y]:x\le F(x)\}.
\end{equation}
Applying these relations one can show that the fixed points $y_*$ and $y^*$ of $F$ are increasing with respect to $F$. Consequently, by (\ref{E30}) and (\ref{E02}), their components, and hence the smallest and greatest solutions of the system (\ref{E1}) in $[\underline y,\overline y]$, are increasing with respect to $f_i$ and $c_i$, $i=1,\dots,m$.    
\qed\vskip6pt
 
As a special case of Theorem \ref{T21}  we obtain the following
corollary. 
 
\begin{corollary}\label{C4.02} Assume that $f_i:HK_{loc}([a,b))^m\to \mathcal A_D([a,b))$ is increasing and order-bounded
for every $i=1,\dots,m$.
Then the system (\ref{E1}) has in $\mathcal D^{lc}([a,b))^m$ the smallest and greatest solutions, and they are increasing with respect to $f_i$ and $c_i$, $i=1,\dots,m$.
\end{corollary}

{\bf Proof}. Because functions $f_i$ are order-bounded, there exist distributions  $\underline h_i$ and $\overline h_i$ in $\mathcal A_D([a,b))$ such that
$\underline h_i\preceq f_i(x_1,\dots,x_m)\preceq \overline h_i$ for all $(x_1,\dots,x_m)\in HK_{loc}([a,b))^m$ and $i=1,\dots,m$. 
Defining for every $i=1,\dots,m$,
$$
\underline y_i(t):=c_i+\int_a^t \underline h_i, \quad \overline y_i(t):=c_i+\int_a^t \overline h_i, \ t\in[a,b),
$$
it is easy to see that $\underline y=(\underline y_1,\dots,\underline y_m)$ is a lower solution
and $\overline y=(\overline y_1,\dots,\overline y_m)$ is an upper solution of (\ref{E1}) in $\mathcal D^{lc}([a,b))^m$.
Thus (\ref{E1}) has by Theorem \ref{T21} the smallest and greatest solutions in the order interval $[\underline y,\overline y]$ of $\mathcal D^{lc}([a,b))^m$.
If $y=(y_1,\dots,y_m)$ is a solution of (\ref{E1}) in $\mathcal D^{lc}([a,b))^m$, then
$$
 \underline y_i(t):=c_i+\int_a^t \underline h_i\le c_i+\int_a^tf_i(y)\le \overline c_i+\int_a^t\overline h_i=\overline y_i(t), \ t\in[a,b),\ i=1,\dots,m.
 $$
Thus $y$ belongs to the order interval  $[\underline y,\overline y]$ of $\mathcal D^{lc}([a,b))^m$, whence the smallest and greatest solutions of (\ref{E1}) in that order interval are the smallest and greatest solutions of (\ref{E1}) in 
the whole $\mathcal D^{lc}([a,b))^m$. The last conclusion of Theorem \ref{T21} implies that these solutions are  increasing with respect to $f_i$ and $c_i$, $i=1,\dots,m$.
\qed\vskip6pt
  
The following result is a consequence of Corollary  \ref{C4.02}. 

\begin{proposition}\label{P21}
Assume that for all $i=1,\dots,m$ and $(x_1,\dots,x_m)\in HK_{loc}([a,b))^m$ 
each $f_i(x_i,\dots, x_m)$ is the distributional derivative of a function 
\begin{equation}\label{E31}
F_i(x_1,\dots,x_m)(t) = \sum_{j=1}^nH_{ij}(t)\,\kint_a^tg_{ij}(x_1,\dots,x_m)+G_i(t), \quad t\in[a,b),
\end{equation}
where  $G_i$ belong to $\mathcal D^{lc}([a,b))$, $G_i(a)=0$, functions  $H_{ij}$ are bounded and non-negative-valued on $[a,b)$ and left-continuous on $(a,b)$, and  the functions $g_{ij}(x):[a,b)\to\mathbb R$ satisfy the following hypotheses:
\begin{description}
\item[(g$_{ij1}$)] $g_{ij}(x_1,\dots,x_m)$ is locally $HK$ integrable for all $(x_1,\dots,x_m)\in HK_{loc}([a,b))^m$.
\item[(g$_{ij2}$)] There exist locally $HK$ integrable functions $\underline g_{ij},\overline g_{ij}:[a,b)\to \mathbb R$ such that 
\item[]\qquad      
$\kint_a^t\underline g_i\le \kint_a^tg_{ij}(x_1,\dots,x_m)\le \kint_a^tg_{ij}(y_1,\dots,y_m)\le \kint_a^t\overline g_{ij}$, whenever $t\in[a,b)$ and 
\item[]\qquad  $(x_1,\dots,x_m)\le (y_1,\dots,y_m)$ in $HK_{loc}([a,b))^m$. 
\end{description}
Then the system (\ref{E1}) has in $\mathcal D^{lc}([a,b))^m$ the smallest and greatest solutions, and they are increasing with respect to $g_{ij}$ and $c_i$.
\end{proposition}

{\bf Proof}. The hypotheses ensure that (\ref{E31}) defines for every $(x_1,\dots,x_m)\in HK_{loc}([a,b))^m$ and $i=1,\dots,m$ a  function $F_i(x_1,\dots,x_m)\in\mathcal D^{lc}([a,b))$, and that its distributional derivative $f_i(x_1,\dots, x_m)$ is increasing in $(x_1,\dots, x_m)$ and is order-bounded by distributions $\underline h_i$ and $\overline h_i$ whose primitives are
$$ 
\underline y_i(t) = c_i+\sum_{j=1}^nH_{ij}(t)\,\kint_a^t\underline g_{ij}+G_i(t), \ \overline y_i(t) = c_i+\sum_{j=1}^nH_{ij}(t)\,\kint_a^t\overline g_{ij}+G_i(t),\quad t\in[a,b).
$$ 
Thus the conclusions follow from Corollary  \ref{C4.02}. 
\qed\vskip6pt

\begin{remark}\label{R21} 
The smallest elements of the well-ordered chain $C$ determined by (I) are $F^n(\underline y)$, $n\in\mathbb N_0$, as long as $F^n(\underline y)=F(F^{n-1}(\underline y))$  is defined and $F^{n-1}(\underline y)<F^n(\underline y)$, $n\in\mathbb N$. If $F^{n-1}(\underline y)=F^n(\underline y)$ for some $n\in\mathbb N$, there
 is a smallest such $n$, and $y_*=F^{n-1}(\underline y)$ is under the hypotheses of  Theorem \ref{T21} the smallest fixed point of $F$ in $[\underline y,\overline y]$. 
If $y_\omega=\underset{n\in\mathbb N}{\sup}F^n(\underline y)$ is defined in $HK_{loc}([a,b))^m$ and is a strict upper bound of  $\{F^n(\underline y)\}_{n\in\mathbb N}$, then $y_\omega$ is the next element of $C$. If $y_\omega=F(y_\omega)$, the $y_*=y_\omega$, otherwise the next elements of $C$ are of the form $F^n(y_\omega)$, $n\in\mathbb N$, and so on.

The greatest elements of the inversely well-ordered chain $D$ determined by (II)  are $n$-fold iterates $F^n(\overline y)$, as long as they are defined and $F^{n}(\overline y)<F^{n-1}(\overline y)$. If equality holds for some $n\in\mathbb N$, then $y^*=
F^{n-1}(\overline y)$ is the greatest fixed point of $F$ in $[\underline y,\overline y]$.
\end{remark}   
 
\begin{example}\label{Ex31} Determine the smallest and greatest solution of   the system (\ref{E1}), where $m=2$,
$c_i=0$,  $f_i(x_1,x_2)$ are for each $(x_1,x_2)\in  HK_{loc}([0,\infty))^2$ the distributional derivatives of the functions $F_i(x_1,x_2):HK_{loc}([0,\infty))^2\to \mathcal D^{lc}([0,\infty))$, defined by
\begin{equation}\label{E34}
 F_i(x_1,x_2)(t)= H_1(t)\,\kint_0^tg_{i1}(x_1,x_2)+G_i(t), \quad t\in[0,1], \ i=1,2,
\end{equation}
where $H_1$ is the Heaviside step function,  $G_i\in \mathcal D^{lc}([0,\infty))$, $G_i(0)=0$, $g_{i1}(x_1,x_2)(0)=0$, $i=1,2$, and
$$
\begin{cases} 
g_{11}(x_1,x_2)(t)=\arctan\left([10^5\kint_0^1(x_2(t)-G_2(t))\,dt]10^{-4}\right)\left(\frac 1t\cos(\frac 1t)-\sin(\frac 1t)+1\right),\\ 
g_{21}(x_1,x_2)(t)=\tanh\left([3\cdot 10^4\kint_0^1(x_1(t)-G_1(t))\,dt]10^{-4}\right)\left(\frac 1t\sin(\frac 1t)+\cos(\frac 1t)+1\right),
\end{cases}
$$
$[z]$ denoting the greatest integer $\le z\in\mathbb R$.
\end{example}

{\bf Solution:} The validity of the hypotheses (g$_{i11}$) and (g$_{i12}$) is easy to verify. Thus, the system (\ref{E1}) has by Proposition \ref{P21} the smallest and greatest solutions in $\mathcal D^{lc}([0,\infty))^2$. 
To determine these solutions, denote
$$
\begin{cases}
\underline y_1(t):=G_1(t)- 4t(1+\cos(\frac 1t)),\ t\in(0,1], \quad \underline x_1(0)=0,\\
\overline y_1(t):=G_1(t)+ 4t(1+\cos(\frac 1t)),\ t\in(0,1], \quad \overline x_1(0)=0,\\
\underline y_2(t):=G_2(t)- 4t(1-\sin(\frac 1t)),\ t\in(0,1], \quad \underline x_2(0)=0,\\
\overline y_2(t):=G_2(t)+ 4t(1-\sin(\frac 1t)),\ t\in(0,1], \quad \overline x_0(0)=0.
\end{cases}
$$
Calculating the successive approximations
$$
\begin{cases}
(x_{n+1},y_{n+1})=(F_1((x_n,y_n),F_2(x_n,y_n)), \quad (x_0,y_0)=(\underline y_1,\underline y_2) \ \hbox{ and } \\ 
(z_{n+1},w_{n+1})=(F_1(z_n,w_n),F_2(z_n,w_n), \quad (z_0,w_0)=(\overline y_1,\overline y_2),
\end{cases}
$$
we see that  $(x_n,y_n)$ form an increasing  and $(z_n,w_n)$ form a decreasing sequence. Moreover, $(x_{16},y_{16})=(F_1(x_{16},y_{16}),F_2(x_{16},y_{16}))$, and that $(z_{16},w_{16})=(F_1(z_{16},w_{16}),F_2(z_{16},w_{16}))$.
According to (I), (II) and Remark \ref{R21} we then have  
$$
C=\{(x_n,y_n)\}_{n=0}^{16},\ \sup F[C]=(x_{16},y_{16}),\ D=\{(z_n,w_n)\}_{n=0}^{16},\ \hbox { and } \ \inf F[D]=(z_{16},w_{16}).
$$
Thus $y_*=(x_{16},y_{16})$ and $y^*=(z_{16},w_{16})$ are the smallest and greatest solutions of  (\ref{E1}) when $f_i(x_1,x_2)$ are for each $(x_1,x_2)\in HK_{loc}([0,\infty))^2$ the distributional derivatives of the functions $F_i(x_1,x_2):HK_{loc}([0,\infty))^2\to \mathcal D^{lc}([0,\infty))$, defined by (\ref{E34}). The exact formulas of $y_*$ and $y^*$, calculated by using simple Maple programs, are
$$\begin{cases}
y_*(t)=\left(G_1(t)-\arctan\left(\frac{5139}{5000}\right)\left(t+t\cos(\frac 1t)\right),G_2(t)-\tanh\left(\frac{12421}{10000}\right)\left(t-t\sin(\frac 1t)\right)\right), \\ 
y^*(t)=\left(G_1(t)+\arctan\left(\frac{2569}{2500}\right)\left(t+t\cos(\frac 1t)\right),G_2(t)+\tanh\left(\frac{12419}{10000}\right)\left(t-t\sin(\frac 1t)\right)\right). \\ 
\end{cases}
$$

\subsection{Uniqueness results}\label{S82}

Denoting $\theta(t)\equiv 0$ and $\mathcal H_+([a,b))=\{\left\lceil x\right\rceil\,|\,x\in\mathcal D^{lc}([a,b))^m\}$,  where $a < b \le\infty$,
\begin{equation}\label{E401}
\left\lceil x\right\rceil = t\mapsto \max\{|x_i(t)|:i=1,\dots,m\}, \ x=(x_1,\dots,x_m)\in\mathcal D^{lc}([a,b))^m,
\end{equation}
we shall prove the following uniqueness Lemma.

\begin{lemma}\label{L11} Given $F:\mathcal D^{lc}([a,b))^m\to \mathcal D^{lc}([a,b))^m$, assume that
\begin{equation}\label{E1.1}
\left\lceil F(y) - F(z)\right\rceil\le G(\left\lceil y - z\right\rceil) \ \hbox{ whenever } \ y,\,z\in \mathcal D^{lc}([a,b))^m, 
\end{equation}
where $G:\mathcal H_+([a,b))\to \mathcal H_+([a,b))$ has the following properties.
\begin{description}
\item[(G)] $G$ is increasing, i.e., $G(u)\le G(v)$ whenever $u\le v$, and for each $u\in \mathcal H_+([a,b))$ there exists a $w_0\in \mathcal H_+([a,b))$, $u\le w_0$, such that $\inf G[W]=\theta$, where
$W$ is the  chain in $\mathcal H_+([a,b))$ that is inversely well-ordered, i.e., every nonempty subset of $W$ has the greatest element, and that satisfies the following condition.
\item[(W)] $\max W=w_0$, and if $u < w_0$, then $u\in W$ if and only if $u=\inf G[\{w\in W:u<w\}]$.
\end{description}   
Then $F$ has at most one fixed point $y$, i.e., $y\in \mathcal D^{lc}([a,b))^m$ and $y=F(y)$.
\end{lemma}

{\bf Proof}. Assume that $y,\,z\in \mathcal D^{lc}([a,b))^m$, $y=F(y)$,  and $z=F(z)$.
Choose
$w_0\in \mathcal D^{lc}([a,b))$ such that $\left\lceil y - z\right\rceil\le w_0$, and that $\inf G[W]=\theta$, where
$W$ is the  chain in $\mathcal H_+([a,b))$ that is inversely well-ordered and satisfies condition (W). By \cite[Proposition 1.2.1]{HL94} $W$ exists and is uniquely determined. If the inequality $\left\lceil y-z\right\rceil\le w$ does not hold for all $w\in W$, there is the greatest
element in $W$, say $u$, for which $\left\lceil y-z\right\rceil\not\le u$. If $w\in W$ and $u < w$, then $\left\lceil y-z\right\rceil\le w$. This inequality, property (\ref{E1.1}), 
monotonicity of $G$, and equations $y=F(y)$ and $z=F(z)$ imply that 
\begin{equation}\label{E1.2}
\left\lceil y-z\right\rceil=\left\lceil F(y)-F(z)\right\rceil\le G(\left\lceil y-z\right\rceil)\le G(w).
\end{equation}
This result holds for all $w\in W$, $u < w$, whence $\left\lceil y-z\right\rceil$ is a lower bound of the set $G[\{w\in W:u<w\}]$. But $u$ is by (W) the greatest
lower bound of $G[\{w\in W:u<w\}]$, so that $\left\lceil y-z\right\rceil\le u$; a contradiction.

By the above proof $\left\lceil y-z\right\rceil\le w$, and hence (\ref{E1.2}) holds for every $w\in W$, whence
$\left\lceil y-z\right\rceil$ is a lower bound of $G[W]$ in $\mathcal H_+([a,b))$. Thus  $\theta\le \left\lceil y-z\right\rceil\le\inf G[W]=\theta$, i.e., $y=z$.
\qed
 
\begin{remark}\label{R00}
Note that $\mathcal H_+([a,b))$ is not a subset of $\mathcal D^{lc}([a,b))^m$ because the functions of $\mathcal D^{lc}([a,b))$
are not absolutely integrable.

The first elements of the chain $W$ satisfying (W) are $n$-fold iterates $w_n=G^n(\max W)$, $n\in\mathbb N_0$. If $w_\omega=\inf\{w_n\}_{n=1}^\infty$ exists and is a strict lower bound of $\{w_n\}_{n=1}^\infty$, then $w_\omega$ is the next element of $W$. The next possible elements of $W$ are of the form $G^n(w_\omega)$, $n\in\mathbb N$, and so on.
\end{remark}

\begin{example}\label{Ex01} Given $H\in \mathcal D^{lc}([0,1]$, define $F:\mathcal D^{lc}([0,\infty))\to \mathcal D^{lc}([0,\infty))$ by
\begin{equation}\label{E001}
 F(x)(t)= \begin{cases}H(t)+ t(x(t)-H(1)), & 0\le t\le 1,\\
                 x(i)+i+(t-i)(x(t)-x(i)-i), &i<t\le i+1, \ i\in\mathbb N.
 \end{cases}                
\end{equation}

(a) Show that $F$ has  at most one fixed point  $y\in \mathcal D^{lc}([0,\infty))$.
\newline
(b) Determine the fixed point of $F$ when  $H'_-(1)=\underset{t\to 1-}{\lim}\frac{H(1)-H(t)}{1-t}$ exists. 
\end{example}

{\bf Solution:} We will show that the hypotheses of Lemma \ref{L11} hold when  the mapping 
 $G:\mathcal H_+(([0,\infty))\to \mathcal H_+(([0,\infty))$ is defined by
\begin{equation}\label{E003}
 G(u)(t)= \begin{cases} tu(t), & 0\le t\le 1,\\
                 (i+1-t)u(i)+(t-i)u(t), &i<t\le i+1, \ i\in\mathbb N. \end{cases}                
\end{equation}
It is easy to verify that $G$ is increasing, and that (\ref{E1.1}) holds.
Given   $u\in \mathcal H_+([0,\infty))$, define 
$$
w_0(t)=\begin{cases} u(t), \ &0\le t\le 1,\\ \max\{u(t),u(i)\},\ & i<t\le i+1, \ i\in\mathbb N.\end{cases}
$$
Routine calculations and induction imply that 
$$
G^n(w_0)(t)=\begin{cases}  t^{n}u(t), & 0\le t\le 1,\\
                  u(i)+(t-i)^n(w_0(t)-u(i)), & i<t\le i+1, \ i\in\mathbb N,
 \end{cases}  \quad n\in\mathbb N.              
$$
 Thus
$$
\lim_{n\to\infty}G^n(w_0)(t)=\begin{cases} 0, & 0\le t < 1,\\ u(i), & i<t\le i+1, \ i\in\mathbb N. 
\end{cases}                  
$$
Redefining the limit so that the obtained function is left-continuous  at $t=1$ we get  the infimum of the set 
$\{G^n(w_0)\}$ in $\mathcal H_+([0,\infty))$: 
$$
z_1(t)= \inf_{n\in\mathbb N}\{G^n(w_0)\}(t)=\begin{cases} 0, & 0\le t \le 1,\\ u(i), & i<t\le i+1, \ i\in\mathbb N.
                  \end{cases}                
$$ 
Similar calculations and induction show that for every $i=2,3,\dots$,
$$
G^n(z_{i-1})(t)=\begin{cases} 0, & 0\le t\le {i-1},\\ 
(t-i-1)^{n}u(i-1), & i-1 < t\le i,\\
                  u(j), & j<t\le j+1, \ j=i,i+1,\dots,
 \end{cases}  \quad n\in\mathbb N,              
$$
and
$$
z_{i}(t)= \inf_{n\in\mathbb N}\{G^n(z_{i-1})\}(t)=\begin{cases} 0, & 0\le t \le i,\\ u(j), &  j<t\le j+1, \ j= i,i+1,\dots.
\end{cases}
$$                            
Finally, we obtain, as $i\to \infty$,
$$
z_\infty(t)= 0, \quad 0\le t<\infty.
$$
Consequently,  if $u(i)> 0$, $i=1,\dots,b$, the members of  the inversely well-ordered chain $W$ satisfying (W)  are  
$$
G^n(w_0), n=0,1,\dots, \  G^n(z_i), n=1,2,\dots, \ i=1,2,\dots, z_\infty.
$$
In particular, $\inf W=z_\infty=\theta$. 

The above calculations and monotonicity of $G$ imply that the hypotheses
of Lemma \ref{L11} hold. Thus $F$ has  at most one fixed point in $\mathcal D^{lc}([0,\infty))$.

If $y$ is a fixed point of $F$ in $\mathcal D^{lc}([0,\infty))$, it follows from (\ref{E001}) that
$$
 y(t)= \begin{cases} H(t)+ t(y(t)-f(1)), \ 0\le t\le 1,\\
                 y(i)+i+(t-i)(y(t)-y(i)-i), \ i<t\le i+1, \ i\in\mathbb N.
                 \end{cases}
$$                 
Thus
$$
y(t)= H(t)-t\frac{H(1)-H(t)}{1-t}, \quad 0 \le t < 1.
$$
Assuming that  $H'_-(1)=\underset{t\to 1-}{\lim}\frac{H(1)-H(t)}{1-t}$ exists, then $x$ is left-continuous at $t=1$, and 
$y(1)=H(1)-H'_-(1)$. 
Moreover, 
$$
 y(t)=     \frac{(1-t+i)(y(i)+i)}{1-t+i}=y(i)+i, \ i<t<i+1, \ i\in\mathbb N.
$$          
This result  and left-continuity of $y$ at $t=i$, $i=1,\dots,b$, imply that $y(i+1)=y(i)+i$, $i\in\mathbb N$.
Since $y(1)=H(1)-H'_-(1)$, then  
$$
y(t)=H(1)-H'_-(1)+\frac{i(i+1)}2, \quad i<t\le i+1, \ i\in\mathbb N.
$$
Thus the fixed point of $F$ in $\mathcal D^{lc}([0,\infty))$ is
\begin{equation}\label{E801}
 y(t)= \begin{cases} H(t)- t\frac{H(1)-H(t)}{1-t}, \quad &0\le t< 1,\\
 H(1)-H'_-(1), \quad  &t=1,\\
                 H(1)-H'_-(1)+\frac{i(i+1)}2, \quad  &i<t\le i+1, \ i\in\mathbb N.
                 \end{cases}
\end{equation}                 

{\bf Consequence}: Let $F$ be defined by (\ref{E001}) with
$H(0)=0$, then $F(y)\in\mathcal D^{lc}_0([0,\infty))$, where $y$ is defined by (\ref{E801}). Denoting by $f(x)$ the distributional derivative of $F(x)$, $x\in\mathcal D^{lc}([0,\infty))$, then
$$
y(t)=F(y)(t)=\int_0^tf(y), \quad t\in[0,\infty). 
$$  
Because $y(0)=0$ then $y$ is the solution of the Cauchy problem
\begin{equation}\label{E00}
y'=f(y),\quad y(0)=0.
\end{equation}

\subsection{Existence of minimal and maximal solutions}\label{S83} 

In this section sufficient conditions are introduced for  the existence of local or global minimal and maximal  solutions to the distributional Cauchy system
\begin{equation}\label{E11}    
y_i'=f_i(y_1,\dots,y_m), \quad  y_i(a)=0, \ i=1,\dots,m. 
 \end{equation}
We assume that $I=[a,b]$, $a < b < \infty$, or $I=[a,b)$, $a < b\le\infty$.
The space $L^1(I)$ of Lebesgue integrable functions on $I$, ordered a.e. pointwise and normed by $L^1$-norm: $\|x\|_1=\int_a^b|x(s)|\,ds$, is an ordered Banach space.  
It is easy to verify that the product space $L^1(I)^m=\times_{i=1}^mL^1(I)$,
ordered by componentwise ordering and a norm  $\|x\|=\{\max \|x_i\|_1:i=1,\dots,m\}$, $x=(x_1,\dots,x_m)\in L^1(I)^m$, has the following properties.
\begin{description}
\item [(L0)] Bounded and monotone sequences of $L^1(I)^m$ converge.
\item [(L1)] $x^+=\sup\{(\theta,\dots,\theta),x\}$ exists, and $\|x^+\|\le \|x\|$ for every $x\in L^1(I)^m$.
\end{description}
 Denote  
\begin{equation}\label{E32}
 B(R)=\{x\in L^1(I)^m:\|x\|\le R\}, \ R>0.
\end{equation}
Because of the properties (L0) and (L1) the following result is a consequence of \cite[Theorem 2.44]{CH11}.

\begin{lemma}\label{L31} 
Given a subset $P$ of $L^1(I)^m$, assume that $F:P\to P$ is increasing, and
that $F[P]\subseteq B(R)\subseteq P$ for some $R> 0$.  Then $F$ has
minimal and maximal fixed points.
\end{lemma}

The next result is a special case of Lemma \ref{L31}.

\begin{proposition}\label{P31} Assume that  mappings $f_i:L^1(I)^m\to \mathcal A_L(I)$ are increasing,  and that  for some $R > 0$ the $LL$ primitive integrals $F_i(x_1,\dots,x_m)(t)=\rint_a^tf_i(x_1,\dots,x_m)$, $t\in I$, of $f_i(x_1,\dots,x_m)$,  $i=1,\dots,m$, satisfy the following hypothesis.
\begin{description}
\item[(f0)] 
 $\|F_i(x_1,\dots,x_m)\|_1\le R$ for all $x_i\in L^1(I)$, $\|x_i\|_1\le R$, $i=1,\dots,m$.
\end{description}
Then the system (\ref{E11}) has 
 minimal and maximal solutions in $B(R)\cap \mathcal L^{lc}(I)^m$.
\end{proposition}

{\bf Proof}. By definition, the functions $F_i(x_1,\dots,x_m)$ belong to $\mathcal L^{lc}(I)$ which is a subset of 
$L^1(I)$. Thus  $F=(F_1,\dots,F_m)$ maps $L^1(I)^m$ into $L^1(I)^m$.
The given hypotheses imply that $F$ satisfies the hypotheses of Lemma \ref{L31} when
$P=B(R)$. Thus, by Lemma \ref{L31}, $F$ has in $B(R)$ minimal and maximal fixed points. Their components form minimal and maximal solutions of (\ref{E11})  in $B(R)\cap \mathcal L^{lc}(I)^m$. 
\qed

The following result deals with the existence of minimal and maximal  solutions of the system (\ref{E11}) 
 in the whole $\mathcal L^{lc}(I)^m$.

\begin{theorem}\label{T32}  Assume that  mappings $f_i:L^1(I)^m\to \mathcal A_L(I)$ are increasing, 
and 
that the  integrals $F_i(x)(t)=\rint_a^tf_i(x)$, $t\in I$,
 satisfy the following hypothesis.
\begin{description}
\item[(f1)] 
 $\|F_i(x)\|_1\le Q(\|x\|)$ for all $x\in L^1(I)^m$, where
$Q:\mathbb R_+\to\mathbb R_+$ is increasing, $R=Q(R)$ for some $R>0$, and $r\le Q(r)$ implies $r\le R$.
\end{description}
Then the  Cauchy system (\ref{E11}) has minimal and maximal solutions in $\mathcal L^{lc}(I)^m$.
\end{theorem}

{\bf Proof}. 
The hypothesis (f1) implies that $F=(F_1,\dots,F_m)$ has the following property.
$$
\|F(x)\|_1\le Q(\|x\|)\le Q(R)=R\ \hbox{ for every } \ x\in B(R). 
$$
Thus the hypothesis (f0) holds, whence
 (\ref{E11}) has the by Proposition \ref{P31} minimal and maximal solutions  in $B(R)\cap \mathcal L^{lc}(I)^m$, and they are increasing with respect to $f_i$.

If $y=(y_1,\dots,y_m)\in \mathcal L^{lc}(I)^m$ is a solution of (\ref{E11}), then $y$ is a fixed point of $F$. The hypothesis (f1) with $r=\|y\|$ implies that
$$
\|y\|=\|F(y)\|\le Q(\|y\|)\le Q(R)=R. 
$$      
Thus $y\in B(R)$, whence all the solutions of (\ref{E11}) are in $B(R)\cap \mathcal L^{lc}(I)^m$.

The assertion follows from the above results.
\qed

\begin{definition}\label{D30} Given a nonempty set $X$ and a normed space $Y=(Y,\|\cdot\|)$, we say that a mapping 
$f:X\to Y$ is norm bounded if $\sup\{\|f(x)\|:x\in X\}< \infty$.
\end{definition} 
It follows from Theorem \ref{T61} that  $\mathcal A_L(I)$, equipped with the 1-norm:
\begin{equation}\label{E33}
\|g\|_1= \smallint_a^b|G(t)|\,dt\ \hbox{ where } \ G(t)= \int_a^tg, \ g\in \mathcal A_L(I),\ t\in I,
\end{equation}
and the partial order $\preceq$ defined by (\ref{E30}), is a normed Riesz space.

The following result is a consequence of Theorem \ref{T32}.

\begin{corollary}\label{C32} Assume that mappings $f_i:L^1(I)^m\to \mathcal A_L(I)$, $i=1,\dots,m$, are increasing and norm bounded. Then the distributional Cauchy system (\ref{E11}) has minimal and maximal solutions in $\mathcal L^{lc}(I)^m$.
\end{corollary}

{\bf Proof}. Equivalent to the norm boundedness of mappings $f_i$ is  that equations
$$
F_i(x_1,\dots,x_m)(t)=\rint_a^tf_i(x_1,\dots,x_m), \quad t\in I, \quad i=1,\dots,m,
$$  
define norm-bounded mappings $F_i:L^1(I)^m\to \mathcal L^{lc}(I)$, $i=1,\dots,m$.
Thus there exist $R_i>0$ such that 
$$
\|F_i(x)\|_1 \le R_i \ \hbox{ for all } \ x\in L^1(I)^m, \ i=1,\dots, m.
$$
These inequalities imply that the hypothesis (f1) is valid when
$$
Q(r)\equiv R:=\max\{R_i|i=1,\dots,m\}.
$$
Because the mappings $f_i$ are also increasing, the conclusion follows from Theorem \ref{T32}.
\qed

\begin{remark}\label{R05} It follows from \cite[Theorem 2.44]{CH11} that the system (\ref{E11}) has under the hypotheses of Proposition \ref{P31} or Theorem \ref{T32} the
smallest solution $y_*=(y_{*1},\dots,y_{*m})$ and the greatest solution $y^*=(y_1^*,\dots,y_m^*)$ in the order
interval $[\underline y,\overline y]$ of $B(R)$,  where $\underline y=(\underline y_1,\dots,\underline y_m)$
is the greatest solution of the system
$$
y_i(t)=\min\{0,\rint_a^tf_i(y_1,\dots,y_m)\}, \quad t\in I, \ i=1,\dots,m,
$$ 
and  $\overline y=(\overline y_1,\dots,\overline y_m)$
is the smallest solution of 
$$
y_i(t)=\max\{0,\rint_a^tf_i(y_1,\dots,y_m)\}, \quad t\in I, \ i=1,\dots,m.
$$
Moreover, $y^*$, $y_*$,  $\underline y$ and  $\overline y$ are all
increasing with respect to $f_i$,  $i=1,\dots,m$.

Existence of continuous solutions of distributional Cauchy problems is studied in \cite{SH110}. 
\end{remark}


\section{Higher order differential equations}\label{S7}
\setcounter{equation}{0}

In this section we will study the following $m$th order  
order distributional Cauchy problem
\begin{equation}\label{E431.01}
\left\{
\begin{aligned}
y^{(m)}&=g(y,y',\dots,y^{(m-1)}),\\
y(0)& = c_{1}, \ y'(0) = c_{2}, \ \dots, \ y^{(m-1)}(0) = c_{m}.
\end{aligned}
\right.
\end{equation}
Existence and comparison
results for the smallest and greatest solutions of problem (\ref{E431.01}) are then proved by using results of subsection \ref{S81}. 
\begin{definition}\label{D431.01}
We say that $y:[0,b)\to \mathbb R$, $b\in (0,\infty]$, is a {\em  solution} of (\ref{E431.01}) if
$y^{(m-1)}\in \mathcal D^{lc}[0,b)$, if $g(y,y',\dots,y^{(m-1)})\in \mathcal A_D[0,b))$,
and if (\ref{E431.01}) holds.
\end{definition}

The Cauchy problem (\ref{E431.01}) can be transformed into a system of first order Cauchy problems
as follows:

\begin{lemma}\label{L431.01} If $g:\mathcal D^{lc}[0,b))^m\to \mathcal A_D[0,b))$, then
$y$ is a  solution of the Cauchy problem
(\ref{E431.01})  if and only if $(y_1,\dots,y_m)=(y,y',\dots,y^{(m-1)})$ is a solution of
the following Cauchy system:
\begin{equation}\label{E431.00}
y_i' =y_{i+1}, \ i=1,\dots,m-1,\ y_m'=g(y_1,\dots,y_m), \  y_i(0) = c_i, \ i=1,\dots,m.
\end{equation}
\end{lemma}

We will present conditions under which problem
(\ref{E431.01}) has the smallest and greatest solutions in the set 
$$
S_D=\{y:[0,b)\to \mathbb R:y^{(m-1)}\in \mathcal D^{lc}[0,b))\}
$$
 or in its order interval when $S_D$ is ordered by
$$
y\le z \ \hbox{iff $y(t)\le z(t)$ and $y^{(i)}(t)\le z^{(i)}(t)$ for all $t\in [0,b)$ and}\ i=1,\dots,m-1.
$$ 
We study also dependence of  solutions of (\ref{E431.01}) on the
functions $g$ and on the initial values $c_i$, $i=1,\dots,m$.
\begin{definition}\label{D433.1}
We say that a function $y\in S_D$ is a {subsolution} of problem
(\ref{E431.01}) if
\begin{equation}\label{E433.02}
y^{(m)}\preceq  g(y,y',y^{(m-1)}),\ 
y(0) \le  c_1, \ y'(0) \le c_{2},  \dots, y^{(m-1)}(0)\le  c_{m}.
\end{equation}
If reversed inequalities hold in (\ref{E433.02}), we say that $y$ is
a {supersolution} of  (\ref{E431.01}). If equalities hold in
(\ref{E433.02}), then $y$ is called a {solution} of
(\ref{E431.01}).
\end{definition}
As a consequence of Theorem \ref{T21} we obtain an existence
comparison theorem for  solutions of problem (\ref{E431.01}).

\begin{theorem}\label{T433.01} Assume that $g:D([0,b))^m\to \mathcal A_D[0,b))$ is increasing, and that
(\ref{E431.01}) has a
subsolution $\underline y\in S_D$ and a supersolution  $\overline y\in
S_D$, and $\underline y\le\overline y,\, \underline y'\le\overline y',
\dots,\underline y^{(m-1)}\le\overline y^{(m-1)}$. Then
the Cauchy problem (\ref{E431.01}) has the smallest and greatest solutions in
the order interval $[\underline y,\overline y]$ of $S_D$, and
they are increasing with respect to $g$ and $c_i$, $i=1,\dots,m$.
\end{theorem}

{\bf Proof}. The function  $\underline z=(\underline y_1,\dots,\underline y_m)=(\underline y ,\underline y',\dots,\underline y^{(m-1)})$ is a subsolution, and the function $\overline z=(\overline y_1,\dots,\overline y_m)=(\overline y ,\overline y',\dots,\overline y^{(m-1)})$ is  a supersolution of
the  system (\ref{E1}) when the functions $f_i:\mathcal D([0,b))^m\to \mathcal A_D[0,b))$ are defined by 
\begin{equation}\label{E431.03}
f_i(x) =x_{i}, \ i=1,\dots,m-1,\ f_m(x)=g(x), \quad x=(x_1,\dots,x_m)\in\mathcal D^{lc}[0,b))^m, 
\end{equation}
These functions $f_i$ are also increasing because $g$ is. It then follows from Theorem \ref{T21} that
the so obtained system (\ref{E1}) has the smallest solution $(y_{*1},\dots,y_{*m})$ and the greatest solution
$(y_1^*,\dots,y_m^*)$ in the order interval $[\underline z,\overline z]$ of $\mathcal D^{lc}([0,b))^m$. Moreover, they are
increasing with respect to $f_i$ and $c_i$, $i=1,\dots,m$.
This result implies by Lemma \ref{L431.01} that $y_{1*}$ and $y_1^*$ are the smallest and the greatest solutions of
the Cauchy problem (\ref{E431.01}) in the order interval $[\underline y,\overline y]$ of $S_D$, and they are increasing with respect to $g$ and $c_i$, $i=1,\dots,m$.
\qed\vskip6pt

Next we  consider the existence of the smallest and greatest
solutions of the Cauchy problem (\ref{E431.01}) in the whole $S_D$.
As a special case of Theorem \ref{T433.01} we obtain the following result.

\begin{corollary}\label{C433.001} The Cauchy problem (\ref{E431.01}) has the smallest and greatest solutions in $S_D$, and they are increasing with respect to $g$ and $c_i$, $i=1,\dots,m$, if   $g:D([0,b))^m\to \mathcal A_D([0,b))$ is increasing and order-bounded.
\end{corollary}

{\bf Proof}. Because function $g$ is order-bounded, there exist distributions  $\underline h$ and $\overline h$ in $\mathcal A_D[0,b))$ such that
$\underline h\preceq g(x_1,\dots,x_m)\preceq \overline h$ for all $(x_1,\dots,x_m)\in D([0,b))^m$. 
Defining,
$$
\begin{cases}
\underline y_m(t):=c_m+\int_0^t \underline h, \quad \overline y_m(t):=c_m+\int_0^t \overline h_i, \ t\in[0,b),\\
\underline y_i(t)=c_i+\int_0^t\underline y_{i+1}(s)\,ds, \quad \overline y_i(t)=c_i+\int_0^t\overline y_{i+1}(s)\,ds,
i=1,\dots,m-1, \ t\in[0,b),
\end{cases}
$$
it is easy to see that $\underline y_1$ is a lower solution
and $\overline y_1$ is an upper solution of (\ref{E431.01}) in $\mathcal D^{lc}([0,b))$.
Thus (\ref{E431.01}) has by Theorem \ref{T433.01} the smallest and greatest solutions in the order interval $[\underline y_1,\overline y_1]$ of $S_D$.
If $y$ is a solution of (\ref{E431.01}) in $S_D$, then
$$
 \underline y_1^{(m-1)}(t)=c_m+\int_0^t \underline h\le c_m+\int_0^tg(y,y',\dots,y^{(m-1)})\le  c_m+\int_0^t\overline h=\overline y_1^{(m-1)}(t), \ t\in[0,b).
 $$
This result and the definitions of $\underline y_1$ and $\overline y_1$ can be used to show that 
$y$ belongs to the order interval  $[\underline y_1,\overline y_1]$ of $S_D$. Thus the smallest and greatest solutions of (\ref{E431.01}) in that order interval are the smallest and greatest solutions of (\ref{E431.01}) in 
the whole $S_D$. The last conclusion of Theorem \ref{T433.01} implies that these solutions are  increasing with respect to $g$ and $c_i$, $i=1,\dots,m$.
\qed

\end{document}